\documentclass[twoside, 12pt]{paper}
\usepackage{a4wide, amsthm, amssymb, amscd, bookman, euler, eucal}

\DeclareSymbolFont{rsfs}{OMS}{rsfs}{m}{n}
\DeclareSymbolFontAlphabet{\mathscr}{rsfs}
\DeclareSymbolFont{bbold}{U}{dsrom}{m}{n}
\DeclareSymbolFontAlphabet{\mathbb}{bbold}

\renewcommand{\rho}{\varrho}

\renewcommand{\Bbb}{\mathbb}
\renewcommand{\frak}{\mathfrak}
\newcommand{\catqot}{/\hskip-3pt/}
\newcommand{\C}{{\Bbb C}}

\newcommand{\E}{{\cal E}}

\newcommand{\F}{{\cal F}}

\newcommand{\GL}{\mathop{\rm GL}}

\newcommand{\Hom}{\mathop{\rm Hom}}

\newcommand{\id}{\mathop{\rm id}}

\renewcommand{\L}{{\cal L}}

\newcommand{\n}{{\cal N}}
\renewcommand{\O}{{\cal O}}

\renewcommand{\P}{{\Bbb P}}

\newcommand{\Q}{{\Bbb Q}}
\newcommand{\R}{{\Bbb R}}

\newcommand{\SL}{\mathop{\rm SL}}

\newcommand{\Spec}{\mathop{\rm Spec}}

\newcommand{\Sym}{\mathop{\rm Sym}}

\newcommand{\Z}{{\Bbb Z}}

\newcommand{\la}{\lambda}
\newcommand{\lra}{\longrightarrow}
\newcommand{\ra}{\rightarrow}

\newcommand{\lma}{\longmapsto}
\newcommand{\p}{\prime}
\newcommand{\q}{\quad}

\renewcommand{\phi}{\varphi}
\newcommand{\rk}{\mathop{\rm rk}}
\newcommand{\eps}{\varepsilon}

\renewcommand{\theta}{\vartheta}
\newcommand{\ul}{\underline}
\newcommand{\ol}{\overline}

\theoremstyle{plain}
\newtheorem{Thm}{Theorem}[subsection]
\newtheorem*{Thm*}{Main Theorem}
\newtheorem{Cor}[Thm]{Corollary}
\newtheorem{Prop}[Thm]{Proposition}
\newtheorem*{Prop*}{Proposition}
\newtheorem{Lem}[Thm]{Lemma}

\newtheorem*{Claim}{Claim}

\theoremstyle{definition}

\theoremstyle{remark}

\newtheorem{Rem}[Thm]{Remark}
\newtheorem{Ex}[Thm]{Example}
\newtheorem*{Rem*}{Remark}

\begin{document}

\pagestyle{myheadings}
\markboth{\rm Alexander H.W.\ Schmitt}{\rm Global Boundedness}

\title{Global Boundedness for Decorated Sheaves\footnote{To appear in Int.\ Math.\ Res.\ Not.}}
\author{Alexander H.W.\ Schmitt
\institution{\sf Universit\"at Duisburg-Essen\\\sf FB6 Mathematik \& Informatik\\\sf D-45117 Essen
\\\sf Germany\\\sf
E-Mail: \tt alexander.schmitt@uni-essen.de\rm}}
\date{}
\maketitle
\begin{abstract}
An important classification problem in Algebraic Geometry deals with pairs $(\E,\phi)$, consisting of
a torsion free sheaf $\E$ and a non-trivial homomorphism $\phi\colon (\E^{\otimes a})^{\oplus b}
\lra\allowbreak  \det(\E)^{\otimes c}\otimes \L$ on a polarized  complex projective
manifold $(X,\O_X(1))$, the input data $a$, $b$, $c$, $\L$ as well as the Hilbert
polynomial of $\E$ being fixed. The solution to the classification problem
consists of a family of moduli spaces 
${\cal M}^\delta:={\cal M}^{\delta-\rm ss}_{a/b/c/L/P}$ for the $\delta$-semistable
objects, where $\delta\in\Q[x]$ can be any positive polynomial of degree at most $\dim X-1$.
In this note we show that there are only finitely many distinct moduli spaces among the ${\cal M}^\delta$
and that they sit in a chain of ``GIT-flips". This property has been known and proved by ad hoc arguments
in several special cases. In our paper, we apply refined information on the instability flag to solve
this problem. This strategy is inspired by the fundamental paper of Ramanan and Ramanathan on the
instability flag.
\end{abstract}
\section{Introduction}
In this paper, we continue the study of the important moduli problem of decorated sheaves \cite{Schmitt0},
\cite{GS}.
The input data for this moduli problem are a polarized complex projective manifold $(X,\O_X(1))$,
non-negative integers $a$, $b$, $c$, and a line bundle $\L$ on $X$. The objects of study are pairs
$(\E,\phi)$ which consist of a torsion free coherent sheaf $\E$ on $X$ and a non trivial homomorphism
$\phi\colon (\E^{\otimes a})^{\oplus b}\lra \det(\E)^{\otimes c}\otimes \L$. We refer to such a pair
$(\E,\phi)$ as a \it torsion free sheaf with a decoration of type $(a,b,c;\L)$ \rm or simply as a \it decorated
torsion free sheaf\rm, if the input data are understood. Two torsion free sheaves $(\E_1,\phi_1)$ and
$(\E_2,\phi_2)$ with a decoration of type $(a,b,c;\L)$ are called \it equivalent\rm, if there are an isomorphism
$\phi\colon\E_1\lra\E_2$ and a number $z\in\C^*$, such that
$$
\phi_2=\bigl(\det(\psi)^{\otimes c}\otimes 
(z\cdot {\id}_\L)\bigr)\circ\phi_1\circ \bigl((\psi^{\otimes a})^{\oplus b}\bigr)^{-1}.
$$
The problem to classify decorated sheaves up to equivalence is related to many classification problems
in Algebraic Geometry. We refer the reader to \cite{BGP} and \cite{Schmitt0} for surveys.
For solving this classification problem in the traditional sense by establishing the existence of a
coarse moduli scheme, one needs to introduce a notion of semistability. One tests the (semi)stability of
a decorated sheaf $(\E,\phi)$ against \it weighted filtrations of $\E$\rm, i.e., against pairs $(\E^\bullet,
\ul{\alpha})$ which consist of a filtration
$$
\E^\bullet:\q 0\subsetneq \E_1\subsetneq\cdots\subsetneq\E_s\subsetneq \E
$$
of $\E$ by saturated subsheaves and a tuple
$$
\ul{\alpha}=(\alpha_1,...,\alpha_s)
$$
of positive rational numbers.
Such a weighted filtration defines the polynomial 
$$
M(\E^\bullet,\ul{\alpha}):=\sum_{i=1}^s\alpha_i\bigl(P(\E)\cdot\rk(\E_i)-P(\E_i)\cdot\rk(\E)\bigr)
$$
with $\Q$-coefficients and the rational number 
$$
\mu(\E^\bullet,\ul{\alpha};\phi)
 :=-\min\Bigl\{\, \gamma_{j_1}+\cdots+\gamma_{j_a}\,|\, (j_1,...,j_a)\in\{\,1,...,s+1\,\}^{\times a}:
\phi_{|(\E_{j_1}\otimes\cdots\otimes\E_{j_a})^{\oplus b}}\not\equiv 0\,\Bigr\}
$$
using $\E_{s+1}=\E$ and
$$
\ul{\gamma}=(\gamma_1,...,\gamma_{r}):=
\sum_{i=1}^s\alpha_i\cdot\bigl(\underbrace{\rk(\E_i)-\rk(\E),...,\rk(\E_i)-\rk(\E)}_{\rk(\E_i)\times},
\underbrace{\rk(\E_i),...,\rk(\E_i)}_{(\rk(\E)-\rk(\E_i))\times}\bigr).
$$
Finally, fix a positive polynomial $\delta\in\Q[x]$ of degree at most $\dim X-1$. Then, a torsion free sheaf
with a decoration of type $(a,b,c;\L)$ is said to be \it $\delta$-(semi)stable\rm, if for every weighted
filtration $(\E^\bullet,\ul{\alpha})$ of $\E$, the inequality
$$
M(\E^\bullet,\ul{\alpha})+\delta\cdot \mu(\E^\bullet, \ul{\alpha};\phi)(\succeq)0
$$
is satisfied. Here, the notation ``$(\succeq)$" means that ``$\succ$" is to be used in the definition
of ``stable" and ``$\succeq$" in the definition of ``semistable", and ``$\succeq$" and ``$\succ$" refer
to the lexicographic ordering of polynomials.
\par
The classification problem for $\delta$-semistable decorated sheaves $(\E,\phi)$
where the Hilbert polynomial of $\E$ is required to be a fixed polynomial $P$ is solved abstractly
by a projective moduli space ${\cal M}^\delta:={\cal M}^{\delta\rm-ss}_{a/b/c/\L/P}$. The existence of
${\cal M}^\delta$ was established by the author over curves \cite{Schmitt0} and by G\'omez and Sols
over manifolds of arbitrary dimension \cite{GS}. 
For concrete applications, it is important to know that there are only finitely many
different moduli spaces occurring among the ${\cal M}^\delta$ and to relate the different moduli
spaces as explicitly as possible (chain of ``GIT-flips"). The first and most striking study in this
direction is Thaddeus's proof of the Verlinde formula  which deals with the case
$a=b=1$ and $c=0$ \cite{Th}.
In this paper, we will prove that the basic features observed by Thaddeus remain true in the general
case as well. This answers the problem raised by the author in \cite{Schmitt0}.
To be precise, we have the following result.
\begin{Thm*}
Fix the input data $a$, $b$, $c$, and $\L$ as well as the Hilbert polynomial $P$.
Then, there is a finite set $\{\,\widehat{\delta}_1,...,\widehat{\delta}_m\,\}$ of rational polynomials
$0=:\widehat{\delta}_0\prec\widehat{\delta}_1\prec\cdots\prec\widehat{\delta}_m\prec
\widehat{\delta}_{m+1}:=\infty\cdot x^{\dim X-1}$ of degree
at most $\dim X-1$, such that, for every torsion free sheaf $(\E,\phi)$ with a decoration of type $(a,b,c;\L)$
and with $P(\E)=P$, the following properties hold true:
\par
{\rm i)} Suppose there is an index $i\in \{\,0,...,m\,\}$ with 
$\widehat{\delta}_i\prec\delta_1\prec\delta_2\prec\widehat{\delta}_{i+1}$.
Then, $(\E,\phi)$ is $\delta_1$-(semi)stable
if and only if it is $\delta_2$-(semi)stable. In particular, there is a canonical isomorphism
$$
{\cal M}^{\delta_1}\cong {\cal M}^{\delta_2}.
$$
\par
{\rm ii)} Assume $\widehat{\delta}_i\prec\delta\prec\widehat{\delta}_{i+1}$
for some index $i\in \{\, 1,...,m-1\,\}$. If $(\E,\phi)$ is $\delta$-semistable, then
$(\E,\phi)$ is also $\widehat{\delta}_i$- and $\widehat{\delta}_{i+1}$-semistable, so that
there are canonical morphisms
$$
{\cal M}^{\delta}\lra {\cal M}^{\widehat{\delta}_i}\hbox{\q and\q}
{\cal M}^{\delta}\lra {\cal M}^{\widehat{\delta}_{i+1}}
$$
Conversely, if $(\E,\phi)$ is $\widehat{\delta}_i$-
or $\widehat{\delta}_{i+1}$-stable, then
$(\E,\phi)$ is also $\delta$-stable. 
\par
{\rm iii)} Suppose $\delta\succ\widehat{\delta}_{m}$. If $(\E,\phi)$ is $\delta$-semistable, it 
is also $\widehat{\delta}_m$-semistable, so that
there is a natural morphism
$$
{\cal M}^{\delta}\lra {\cal M}^{\widehat{\delta}_m}.
$$
Conversely, if $(\E,\phi)$ is $\widehat{\delta}_m$-stable, then
$(\E,\phi)$ is also $\delta$-stable. 
\par
{\rm iv)} Suppose $0\prec\delta\prec\widehat{\delta}_1$. If $(\E,\phi)$ is $\delta$-semistable,
then $\E$ is a semistable sheaf. 
Letting $\widehat{\cal M}_0$ be the moduli space of semistable sheaves with Hilbert polynomial $P$, we find
a canonical morphism
$$
{\cal M}^\delta \lra \widehat{\cal M}_0.
$$
If $\E$ is a stable sheaf, then $(\E,\phi)$ is $\delta$-stable.
\end{Thm*}
\begin{Rem*}
Observe that, for any positive polynomial $\delta\in\Q[x]\setminus\{\, \widehat{\delta}_1,...,
\widehat{\delta}_m\,\}$ of degree at most $\dim X-1$,
one of the assumptions in ii), iii), or iv) must be satisfied, because the lexicographic
ordering ``$\prec$" defines a total ordering on $\Q[x]$.
\end{Rem*}
We set $\widehat{\cal M}_i:={\cal M}^{\widehat{\delta}_i}$, $i=1,...,m$, ${\cal M}_i:=
{\cal M}^\delta$ for some $\delta$ with $\widehat{\delta}_{i-1}\prec\delta\prec\widehat{\delta}_i$,
$i=1,...,m$, and ${\cal M}_\infty:=
{\cal M}^\delta$ for some $\delta$ with $\delta\succ\widehat{\delta}_m$. Our theorem is then summarized
by the following picture
$$
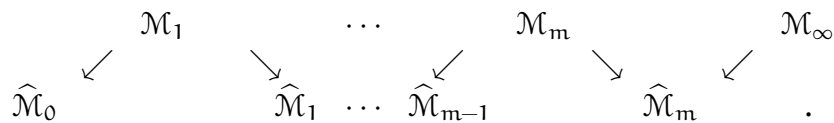
\begin{array}{ccccccccc}
          &         &{\cal M}_1&                 &\cdots
          &         &{\cal M}_{m}\qquad\qquad& &{\cal M}_\infty
\\
          &\swarrow &                      & \searrow &
         &\swarrow                 &\searrow&\swarrow
\\
\widehat{\cal M}_0        &         &                     &\qquad \widehat{\cal M}_1
&\cdots&\widehat{\cal M}_{m-1}& \qquad\qquad\widehat{\cal M}_m&&.
\end{array}
$$
Note that the morphisms occurring, such as ${\cal M}_i\lra \widehat{\cal M}_i$ are induced by taking
quotients of the same projective variety but with respect to different linearizations. 
This fact makes one expect that one may analyze them quite explicitly in a given example 
(see \cite{Th}, \cite{Th2}, \cite{DH}).
In the course of proving the above theorem, we will also work out an explicit description of
the condition of $\delta$-semistability for $\delta\succ\widehat{\delta}_m$.
\par
We also point out that our result does not follow from the ``finiteness of GIT quotients".
Recall that, given an action $\alpha\colon G\times Q\lra Q$ of the reductive group $G$ on the scheme $Q$
and a linearization $\sigma\colon G\times L\lra L$ of $\alpha $ in the line bundle $L$ over $Q$, 
we have the open
subset $Q_\sigma^{\rm ss}$ of points which are semistable w.r.t.\ linearization $\sigma$ and the
GIT quotient $M_\sigma=Q_\sigma^{\rm ss}\catqot G$.
If $Q$ is projective, then one knows that there are only finitely many open subsets $U\subset Q$
of the form $Q_\sigma^{\rm ss}$ and thus only finitely many non-isomorphic quotients
$M_\sigma$, $\sigma$ a linearization of $\alpha$. In the case that $X$ is normal,
this was proved by Dolgachev and Hu \cite{DH}  and, independently and in a more general context,
by Bia\l ynicki-Birula \cite{BB}. In fact, applying Bia\l ynicki-Birula's strategy in the classical 
GIT setting yields a fairly elementary proof for the above fact \cite{FinEasy}.
Of course, for any given polynomial $\delta$, we may find $Q_\delta$, $\alpha_\delta$, and $\sigma_\delta$, 
such that ${\cal M}^\delta=M_{\sigma_\delta}$. But, as can be seen from the constructions
in \cite{GS} and \cite{Schmitt0}, the space $Q^\delta$ depends very much on the specific polynomial $\delta$,
so that the finiteness of GIT quotients may simply not be applied. Only a posteriori with the Main Theorem
at hand may we conclude that all the ${\cal M}^\delta$ may --- up to canonical isomorphy --- be constructed
as quotients of the same projective variety $Q$, because we need to consider only finitely many $\delta$'s. 
That fact, as remarked before, may be useful for studying explicit questions on the moduli spaces.
Our theorem should be considered the analog of the finiteness of GIT quotients within the theory of
decorated sheaves.
\par
Our method of proof relies on the theory of the instability flag developed by Kempf \cite{Kempf} 
and applied by Ramanan and Ramanathan \cite{RamRam}. We establish a certain property of the
instability flag which enables us to modify the strategy of Ramanan and Ramanathan according to our
needs. Another ingredient is the author's analysis of semistability in \cite{Schmitt0}.
Some of the material presented here is also contained in the author's former paper \cite{Schmitt}.
Since we need slightly different notation or an additional argument at some place, we have
repeated several sections in order to make this paper more readable and self-contained.
\section*{Conventions}
We work over the field of complex numbers.
A \it scheme \rm will be a scheme of finite type over $\C$.
For a vector bundle $\E$ over a scheme $X$, we set $\P(\E):={\rm Proj}(\Sym^*(\E))$, i.e., $\P(\E)$ is the projective
bundle of hyperplanes in the fibres of $\E$. An open subset $U\subset X$ is said to be \it big\rm, if
${\rm codim}_X(X\setminus U)\ge 2$.
\par
The degree $\deg(\E)$, the slope $\mu(\E):=\deg(\E)/\rk\E$, and the Hilbert polynomial $P(\E)$
of an $\O_X$-module $\E$ are computed w.r.t.\ the given ample line bundle $\O_X(1)$.
\par
If $\E$ is a torsion free coherent $\O_X$-module, we set
$$
\mu_{\rm max}(\E):=\max\bigl\{\, \mu(\F)\,|\, 0\subsetneq\F\subsetneq\E\,\bigr\}.
$$
If $(\E^\bullet, \ul{\alpha})$ is a weighted filtration of the torsion free sheaf $\E$, then we set
$$
L(\E^\bullet,\ul{\alpha}):=\sum_{i=1}^s \alpha_i \bigl(\deg(\E)\cdot \rk\E_i-\deg(\E_i)\cdot\rk\E\bigr).
$$
This is the coefficient of $x^{\dim X-1}$ in $M(\E^\bullet,\ul{\alpha})$.
\section*{Acknowledgment}
The author acknowledges support by the DFG through a Heisenberg fellowship and through the priority program
``Globale Methoden in der komplexen Geometrie --- Global Methods in Complex Geometry''. 
\par
Most of this paper was conceived and written down during the author's visit to the 
Consejo Superior de Investigaciones Cient\'\i ficas (CSIC) in Madrid
which was funded by the European Differential Geometry Endeavour (EDGE),
EC FP5 contract no.\ HPRN-CT-2000-00101. 
The author wishes to thank O.\ Garc\'\i a-Prada for the invitation and hospitality.
\section{Preliminaries}
\subsection{Geometric Invariant Theory}
\label{GIT}
We use the following convention:
Let $G$ be a complex reductive group which acts on the projective scheme $X$, and suppose this action is linearized
in the ample line bundle $\L$. Given a one parameter subgroup $\la\colon \C^*\lra G$ and a point $x\in X$,
we form $x_\infty:=\lim_{z\ra \infty}\la(z)\cdot x$. Then, $x_\infty$ remains fixed under the $\C^*$-action
induced by $\la$ and the $G$-action, so that $\C^*$ acts on $\L\langle x\rangle$ by a character, say, $z\lma z^\gamma$,
$z\in\C^*$. One sets
$$
\mu_\L(\la,x) :=  -\gamma.
$$ 
\paragraph{One Parameter Subgroups and Parabolic Subgroups. ---}
Let $G$ be a complex reductive group, and $\la\colon \C^*\lra G$ a one parameter subgroup.
Then, we define the parabolic subgroup
\begin{equation}
\label{parabolic}
Q_G(\la) := \Bigl\{\, g\in G\,|\, \lim_{z\rightarrow\infty} \la(z)\cdot g\cdot\la(z)^{-1} \hbox{ exists in $G$}
\,\Bigr\}.
\end{equation}
In fact, any parabolic subgroup of $G$ arises in this way. We refer the reader to the books \cite{Spri} and
\cite{GIT}, Chapter 2.2, for more details. The centralizer $L_G(\la)$ of $\la$ is a \it Levi-component 
of $Q_G(\la)$\rm, i.e., $Q_G(\la)={\cal R}_u(Q_G(\la))\rtimes L_G(\la)$. In this picture, the unipotent radical
of $Q_G(\la)$ is characterized as
$$
{\cal R}_u(Q_G(\la)) = \Bigl\{\, g\in Q_G(\la)\,|\, 
\lim_{z\rightarrow\infty} \la(z)\cdot g\cdot\la(z)^{-1} =e
\,\Bigr\}.
$$
\begin{Rem}
\label{weighted}
i) In the sources quoted above, one takes the limit $z\rightarrow 0$ in order to define a parabolic subgroup
$P_G(\la)$. Thus, we have
\begin{equation}
\label{Sprcomparison}
Q_G(\la)= P_G(-\la).
\end{equation}
\par
ii) Let $G$ be a complex reductive group which acts on the projective scheme $X$, and suppose this action is linearized
in the ample line bundle $\L$. Then, for any point $x\in X$, any one parameter subgroup $\la\colon\C^*\lra G$,
and any $g\in Q_G(\la)$
$$
\mu_\L(\la, x) = \mu_\L(\la,g\cdot x).
$$
This is proved in \cite{GIT}, Chapter 2.2.
\par
iii)
If we are given an injective homomorphism $\iota\colon G\hookrightarrow H$, then we obviously find
$$
Q_H(\la)\cap G= Q_G(\la).
$$
\par
iv) If $G=\GL(V)$, then the group $Q_G(\la)$ is the stabilizer of the flag
$$
V^\bullet\colon\ 0\subsetneq V_1\subsetneq V_2\subsetneq\cdots\subsetneq V_s\subsetneq V 
$$
where $V_i:=\bigoplus_{j=1}^i V^j$, $V^j$ is the eigenspace of the $\C^*$-action coming from $\la$ for
the character $z\lma z^{\gamma_j}$, and $\gamma_1<\cdots<\gamma_{s+1}$ are the different weights occurring.
We also set $\alpha_i:=(\gamma_{i+1}-\gamma_i)/\dim(V)$, $i=1,...,s$. The pair $(V^\bullet,\ul{\alpha})$
is referred to as the \it weighted flag of $\la$\rm. Note, that if $\la^\p$ is conjugate to $\la$, then
$\dim(V_i^\p)=\dim(V_i)$ and $\alpha_i^\p=\alpha_i$, $i=1,...,s$.
\end{Rem}
\paragraph{The Instability Flag. ---}
In this section, ${\Bbb K}$ will be an algebraically closed field of characteristic zero. (Besides for $\C$,
we will need the results also for the algebraic closure of the function field of $X$.)
We start with the group $\GL_n({\Bbb K})$. Let $T$ be the maximal torus of diagonal matrices.
The characters $e_i\colon {\rm diag}(l_1,...,l_n)\lma l_i$, $i=1,...,n$, form a basis for the character
group $X^*(T)$, and 
$$
\begin{array}{rccc}
(.,.)^*\colon & X^*_\R(T)\times X^*_\R(T) &\lra& \R
\\
& \Bigl(\sum_{i=1}^n x_i\cdot e_i, \sum_{i=1}^n y_i\cdot e_i\Bigr) &\lma &\sum_{i=1}^n x_iy_i
\end{array}
$$
defines a scalar product on $X^*_\R(T):=X^*(T)\otimes_\Z\R$ which is invariant under the action of the Weyl group
$W(T):={\cal N}(T)/T$. This yields isomorphisms
$$
X^*_\R(T) \cong  {\Hom}_\R(X^*_\R(T),\R)\cong  X_{*,\R}(T) := X_*(T)\otimes_\Z\R.
$$
For the second identification, we use the duality pairing $\langle.,.\rangle_\R\colon X_{*,\R}(T) \times X^*_\R(T)
\lra \R$ which is the $\R$-linear extension of the canonical pairing 
$\langle.,.\rangle\colon X_{*}(T) \times X^*(T)
\lra \Z$. Since the pairing $(.,.)^*$ is $W(T)$-invariant, the norm $\|.\|_*$ induced on $X_{*,\R}(T)$ extends
to a $\GL_n({\Bbb K})$-invariant norm $\|.\|$ on the set of all one parameter subgroups of $\GL_n({\Bbb K})$ (see \cite{GIT},
Chapter 2.2, Lemma 2.8).
\par
Next, suppose we are given a representation $\kappa\colon \GL_n({\Bbb K})\lra \GL(W)$. This leads to a decomposition
$$
W \cong \bigoplus_{\chi\in X^*(T)} W^\chi
$$
of $W$ into eigenspaces and defines the \it set of weights of $\kappa$ (w.r.t.\ $T$)\rm
$$
{\rm WT}(\kappa,T) := \bigl\{\,\chi\in X^*(T)\,|\, W^\chi\neq\{0\}\,\bigr\},
$$
and, for any $w\in W$, the \it set of weights of $w$ (w.r.t.\ $T$)\rm
$$
{\rm WT}(w,T) := \bigl\{\,\chi\in {\rm WT}(\kappa,T)\,|\, w \hbox{ has a non-trivial component in
}W^\chi\,\bigr\}.
$$
For a one parameter subgroup $\la\in X^*_\R(T)$, we then set
$$
\mu_\kappa(\la,w) := \max\bigl\{\,\langle\la,\chi\rangle_\R\,|\, \chi\in {\rm WT}(w,T)\,\bigr\}.
$$
For any other maximal torus $T^\p\subset G$, we choose an element $g\in G$ with $g\cdot T^\p\cdot g^{-1}=T$, and
set, for $\la\in X^*_\R(T^\p)$,
\begin{equation}
\label{DeFine}
\mu_\kappa(\la,w) :=  \mu_\kappa(g\cdot\la\cdot g^{-1}, g\cdot w).
\end{equation}
\begin{Ex}
i) Let $\P(W^\vee)$ denote the space of lines in $W$. Then, $\kappa$ yields an action of $\GL_n({\Bbb K})$ on 
$\P(W^\vee)$
and a linearization of that action in $\O_{\P(W^\vee)}(1)$. With the former notation, we find
$$
\mu_{\O_{\P(W^\vee)}(1)}(\la, [w]) = \mu_\kappa(\la, w),
$$
for every point $w\in W\setminus\{0\}$ and every one parameter subgroup 
$\la\colon {\Bbb G}_m({\Bbb K})\lra \GL_n({\Bbb K})$.
\par
ii) Our convention is the same as in \cite{Schmitt0}, but differs from the one in
\cite{RamRam}. More precisely, let $\mu^{\rm RR}_\kappa(\la, w)$ be the quantity defined in \cite{RamRam}.
Then,
\begin{equation}
\label{RRcomparison}
\mu_\kappa(\la, w) =  -\mu_\kappa^{\rm RR}(-\la,w).
\end{equation}
\end{Ex}
Now, suppose we are also given a reductive subgroup $G\subset \SL_n({\Bbb K})$. For simplicity, assume that there is a 
maximal torus $T_G$ of $G$ which is contained in  $T$. Otherwise, we may pass to a different maximal torus 
$T^\p$ of $\GL_n({\Bbb K})$.
From $(.,.)^*$ and the dual pairing $(.,.)_*\colon\ X_{*,\R}(T)\times X_{*,\R}(T)\lra \R$, we obtain the
induced pairing $(.,.)_{*,G}\colon\ X_{*,\R}(T_G)\times X_{*,\R}(T_G)\lra \R$. Let $\|.\|_G$ be the restriction of the
norm $\|.\|$ to the one parameter subgroups of $G$. Note that, for $\la\in X_{*,\R}(T_G)$, one
has $\|\la\|_G=\sqrt{(\la,\la)_{*,G}}$.
This last observation implies that $(.,.)_{*,G}$ is invariant under the action of the Weyl group
$W(T_G):=\n_G(T_G)/T_G$. By polarization, this is equivalent to the fact that $\|.\|_G$ restricted
to $X_{*,\R}(T_G)$ is invariant under $W(T_G)$, and this is obvious from the definition.
\begin{Thm}[Kempf]
\label{kempfi}
Suppose $w\in W$ is a $G$-unstable point. Then, the function $\la\lma \nu_\kappa(\la,w):=\mu_\kappa(\la, w)/\|\la\|_G$ 
on the set of all one
parameter subgroups of $G$ attains a minimal value $m_0\in \Q_{<0}$, and there is a unique parabolic subgroup 
$Q(w)\subset G$,
such that $Q(w)=Q_G(\la)$ for every one parameter subgroup $\la\colon {\Bbb G}_m({\Bbb K})\lra G$ with $\nu(\la,w)=m_0$. 
Moreover, if $\la$ and $\la^\p$ are two indivisible one parameter subgroups with 
$\nu(\la,w)=m_0=\nu(\la^\p,w)$, then there
exists a unique element $u\in {\cal R}_u(Q(w))$, such that 
$\la^\p=u\cdot \la\cdot u^{-1}$.
\end{Thm}
\begin{proof} This is Theorem~2.2 in \cite{Kempf}. It is also proved in \cite{RamRam}, Theorem 1.5.
One has to use (\ref{Sprcomparison}) and (\ref{RRcomparison}) to adapt the formulation to the conventions we use.
Since this theorem plays such a crucial r\^ole in our considerations, we briefly remind the reader of the idea of proof.
Recall Equation (\ref{DeFine}) and the fact that $\|\la\|_G=\|g\cdot\la\cdot g^{-1}\|_G$ for all
$\la\colon {\Bbb G}_m({\Bbb K})\lra G$, $g\in G$. First, for an element $g\in G$,
we search for
\begin{equation}
\label{tominimize}
\min\Bigl\{\,\frac{\mu_\kappa(\la,g\cdot w\cdot g^{-1})}{\|\la\|_G}\,\Bigl|\, \la\in X_*(T_G)\,\Bigr\}.
\end{equation}
Write
$$
{\rm WT}(g\cdot w, T_G) := 
= \bigl\{\,\chi^g_1,...,\chi^g_{s(g)}\,\bigr\}.
$$
We obtain the linear forms 
\begin{eqnarray*}
l^g_i\colon\q X_{*,{\Bbb \R}}(T_G) &\lra& {\Bbb \R}
\\
\la &\lma& \langle\la, \chi^g_i\rangle_{\Bbb R},\q i=1,...,s(g),
\end{eqnarray*}
on $X_{*,{\Bbb R}}(T_G)$ which are actually defined over $\Q$.
One has now to study the function 
$$
l^g\colon \la\lma \max_{i=1,...,s(g)} l^g_i(\la)
$$
on the norm-one hypersurface $H$ in $X_{*,{\Bbb R}}(T_G)$ where the assumption is that $l$ possesses
a negative value. One then shows that a function like $l^g$ admits indeed a minimum in a unique point $h\in H$.
Moreover, the fact that the $l^g_i$ are defined over $\Q$ grants that the ray ${\Bbb R}_{>0}\cdot h$
contains rational and integral points. See Lemma 1.1 in \cite{RamRam} for this discussion.
Thus, the expression (\ref{tominimize}) agrees with
$l^g(h)$.
\par
Finally, one remarks that $l^g$ depends only on the set of weights ${\rm WT}(g\cdot w, T_G)$ for which there are only
finitely many possibilities, so that
there is a finite set $\Gamma\subset G$ with
$$
{\rm WT}(g\cdot w, T_G)\in \Bigl\{\,{\rm WT}(\gamma\cdot w, T_G)\,|\,\gamma\in \Gamma\,\Bigr\},\q \hbox{for all $g\in G$}.
$$
Thus, we have to show that 
$$
\min_{\gamma\in \Gamma}\min\left\{\, \frac{l^\gamma(\la)}{\|\la\|_G}\,\Big|\, \la\in X_{*,\R}(T_G)\,\right\}
$$
exists, but this is now clear.
\end{proof}
Let $w$ and $m_0$ be as in the theorem.
We call an indivisible one parameter subgroup $\la\colon{\Bbb G}_m({\Bbb K})\lra G$ with $\nu(\la,w)=m_0$ an \it instability one
parameter subgroup for $w$\rm. Note that, by the theorem, every maximal torus of $Q(w)$ contains a unique instability
one parameter subgroup for $w$.
\begin{Rem}
\label{FiniteInstab}
i) There is also a canonical parabolic subgroup $Q_{\GL_n({\Bbb K})}(w)$ of $\GL_n({\Bbb K})$
with $Q_{\GL_n({\Bbb K})}(w)\cap G=Q(w)$.
Indeed, if $\la$ is any instability subgroup of $w$, then we set $Q_{\GL_n({\Bbb K})}(w):=Q_{\GL_n({\Bbb K})}(\la)$.
This is well-defined because of the last statement in the theorem.
\par
ii) Note that, since ${\rm WT}(g\cdot w, T_G)\subset {\rm WT}(\kappa,T_G)$, and the latter is a finite set,
there are only finitely many possibilities for ${\rm WT}(g\cdot w, T_G)$, so that there are only finitely many
(negative) numbers of the form $m_0$ as $w$ varies over the instable points in $W\setminus\{0\}$ and
$\la$ over the instability one parameter subgroups for $w$. Likewise, by Remark \ref{weighted}, iv), 
the set of data $(\dim V_1,...,\dim V_s; \alpha_1,...,\alpha_s)$ arising from weighted filtrations associated
with instability one 
parameter subgroups of points $w\in W\setminus\{0\}$ is finite. By construction, $\|\la\|_G$ depends only on the
datum $(\dim V_1,...,\dim V_s; \alpha_1,...,\alpha_s)$, whence the set of numbers arising as $\|\la\|_G$ from
an instability one parameter subgroup $\la$ for a point $w\in W\setminus\{0\}$ is finite, too.
\end{Rem}
For every maximal torus $T^\p$ of $\GL_n({\Bbb K})$, the given product on $X^*_\R(T)$ induces the pairing 
$(.,.)^*_{T^\p}\colon X^*_\R(T^\p)\times X^*_\R(T^\p)\lra \R$,
$(\chi,\chi^\p)\lma (\chi(g\cdot .\cdot g^{-1}),\chi^\p(g\cdot .\cdot g^{-1}))^*$, where $g\in\GL_n({\Bbb K})$ is an element,
such that $g\cdot T\cdot g^{-1}=T^\p$. Here, the invariance of $(.,.)^*$ under the Weyl group $W(T)$ implies that 
this product does not depend on the choice of $g$. We set $H_G(w):= Q(w)/{\cal R}_u(Q(w))$, and 
$H_{\GL_n({\Bbb K})}(w):=Q_{\GL_n({\Bbb K})}(w)/{\cal R}_u(Q_{GL(W)}(w))$. Now, $\la$ defines an antidominant
character on $H_{\GL_n({\Bbb K})}(w)$ as follows: Let $\ol{T}$ be a maximal torus of $H_{\GL_n({\Bbb K})}(w)$. Under the isomorphism
$L_{\GL_n({\Bbb K})}(\la)\lra H_{\GL_n({\Bbb K})}(w)$ induced by the quotient morphism $\pi\colon Q_{\GL_n({\Bbb K})}(w)\lra H_{\GL_n({\Bbb K})}(w)$, 
there is a un\-ique maximal torus
$T^\p\subset L_{\GL_n({\Bbb K})}(\la)$ mapping onto $\ol{T}$. Then, as we have explained before, there is a scalar
product $(.,.)_{T^\p}^*\colon X^*_\R(T^\p)\times X^*_\R(T^\p)\lra \R$. This provides us with the unique element 
$l_{T^\p}(\la)$, such that $(l_{T^\p}(\la),\chi)^*_{T^\p}=\langle\la,\chi\rangle$ for all
$\chi\in  X^*_\R(T^\p)$. The computation below 
(Example~\ref{CharComp}) shows 
that $l_{T^\p}(\la)$ is indeed a character of $L_{\GL_n({\Bbb K})}(\la)$ and, thus, of $H_{\GL_n({\Bbb K})}(w)$. Call this character
$\chi_0$. Let $T^{\p\p}$ be any other maximal torus of $Q_{\GL_n({\Bbb K})}(w)$. Then, there is an element
$p\in Q_{\GL_n({\Bbb K})}(w)$ with $p\cdot T^{\p}\cdot p^{-1}=T^{\p\p}$. For all one parameter subgroups 
$\widetilde{\la}\colon {\Bbb G}_m({\Bbb K})\lra T^\p$, 
we have
$$
\langle p\cdot \widetilde{\la}\cdot p^{-1}, \chi_0\rangle
\q=\q
\langle \widetilde{\la},\chi_0\rangle
\q=\q
(\widetilde{\la},\la)^*_{T^\p}
\q=\q
(p\cdot \widetilde{\la}\cdot p^{-1},p\cdot\la\cdot p^{-1})^*_{T^{\p\p}},
$$
so that $p\cdot\la\cdot p^{-1}$ and the maximal torus $\ol{T}^\p:=\pi(T^{\p\p})$ yield indeed the same character
$\chi_0$.
\begin{Ex}
\label{CharComp}
Fix integers $0=:n_0<n_1<\cdots<n_s<n_{s+1}:=n$ and $\gamma_1<\cdots<\gamma_{s+1}$
with $\sum_{i=1}^{s+1} \gamma_{i}(n_{i}-n_{i-1})=0$. This defines a one parameter subgroup $\la\colon
{\Bbb G}_m({\Bbb K})\lra \SL_{n}({\Bbb K})$ via
$$
\la(z)\cdot b_j:= z^{\gamma_i}\cdot b_j,\q j=n_{i-1}+1,...,n_i,\ i=1,...,s+1.
$$
Here, $b_1,...,b_n$ is the standard basis for ${\Bbb K}^n$. Then, $L_{\GL_n({\Bbb K})}(\la)\cong \GL_{n_1}({\Bbb K})\times
\GL_{n_2-n_1}({\Bbb K})\times \cdots\times \GL_{n-n_s}({\Bbb K})$, the latter group being embedded as a group of block diagonal
matrices into $\GL_n({\Bbb K})$. One checks that 
$$
l_T(\la)(m_1,...,m_{s+1})\q=\q\det(m_1)^{\gamma_1}\cdot...\cdot \det(m_{s+1})^{\gamma_{s+1}},\q
\forall\ (m_1,...,m_{s+1})\in L_{\GL_n({\Bbb K})}(\la).
$$
\end{Ex}
Let $w\in W\setminus\{0\}$ be an unstable point, and let $Q(w)\subset G$ be the associated parabolic subgroup.
Moreover, choose an instability one parameter subgroup $\la\colon\ {\Bbb G}_m({\Bbb K})\lra G$ for $w$. This yields, in particular,
a flag $W^\bullet\colon\ 0\subsetneq W_1\subsetneq \cdots\subsetneq W_{t}\subsetneq W$. Next, set
$j_0:=\min\{\, j=1,...,t+1\,|\, w\in W_j\,\}$. Then, $w$ defines a point 
$x_\infty\in \P\bigl( (W_{j_0}/W_{j_0-1})^\vee\bigr)$. Let $m_0\in \Q_{<0}$ be as in Theorem~\ref{kempfi}, 
and $q:=m_0\cdot \|\la\|_G=\mu_\kappa(\la,w)\in\Z_{<0}$. Finally, define $\chi_*:= q\cdot \chi_{0|H_G(w)}$.
\begin{Prop}[Ramanan-Ramanathan]
\label{RamRamI}
The point $x_\infty\in \P\bigl( (W_{j_0}/W_{j_0-1})^\vee\bigr)$ is semistable for the induced $H_G(w)$-action
and its linearization in $\O_{\P( (W_{j_0}/W_{j_0-1})^\vee)}(1)$ twisted by the character $\chi_*$.
\end{Prop}
\begin{proof} This is Proposition 1.12 in \cite{RamRam}. We observe that, by 
(\ref{Sprcomparison}) and (\ref{RRcomparison}),
we have $\chi_*=\chi$ with $\chi$ the character constructed in  \cite{RamRam}.
(Our explicit construction shows that we may take $s=1$ and $r=1$ in the proof of \cite{RamRam}, Proposition 1.12).
Note that Ramanan and Ramanathan
show that $x_\infty=[w_\infty]$, where 
$w_\infty\in W_{j_0}/W_{j_0-1}\otimes {\Bbb K}_{\chi_*^{-1}}$ is a semistable point and ${\Bbb K}_{\chi_*^{-1}}$ 
is the
one dimensional $H_G(w)$-module associated with the character $\chi_*^{-1}$. This gives the claimed linearization.
\end{proof}
\begin{Rem}
\label{FiniteInstabII}
By Remark \ref{FiniteInstab}, ii), it is clear that there are only finitely many possibilities for the
number $q$ introduced above.
\end{Rem}
Finally, we need Kempf's rationality result. For this, let $K$ be a non-algebrai\-cally closed field of characteristic
zero (in our application, this will be the function field of an algebraic variety), $G\lra \Spec(K)$ a $K$-group,
and $W$ a finite dimensional $K$-vector space. Fix an algebraic closure ${\Bbb K}$ of $K$, and set $G_{\Bbb K}
:= G\times_{\Spec(K)}\Spec(\Bbb K)$ and $W_{\Bbb K}:=W\otimes_K{\Bbb K}$. Suppose that we are
given a $K$-rational representation $\kappa\colon\ G_{\Bbb K}\lra \GL(W_{\Bbb K})$.
\begin{Thm}[Kempf]
\label{rationality}
If $T\subset G_{\Bbb K}$ is a maximal torus which is defined over $K$ and $K_T/K$ is a finite extension
of $K$, such that $T_K\times_{\Spec(K)}\Spec(K_T)\cong {\Bbb G}_a(K_T)\times\cdots\times{\Bbb G}_a(K_T)$,
$T_K\subset G$ being the $K$-group with $T_K\times_{\Spec(K)}\Spec(\Bbb K)=T$, then, for a product
$(.,.)^*\colon X_\R^*(T)\times X^*_\R(T)\lra \R$ which is invariant under both the Weyl group $W(T)$ and
the action of the Galois group ${\rm Gal}({\Bbb K}/K)$ via its finite quotient ${\rm Gal}(K_T/K)$, 
the following holds true: 
If $w\in W_{\Bbb K}$ is an unstable $K$-rational point, then the
parabolic subgroup $Q_{G_{\Bbb K}}(w)$, associated to $w$ by means of the norm $\|.\|_{G_{\Bbb K}}$ 
on the one parameter subgroups of $G_{\Bbb K}$ which is induced by $(.,.)^*$, is defined over $K$.
\end{Thm}
\begin{proof}
This is part of Theorem 4.2 in \cite{Kempf}.
See also \cite{RamRam} for generalizations.
\end{proof}
\paragraph{The Instability Flag in a Product. ---}
Here, we begin with two representations $\kappa_1\colon G\lra W_1$ and $\kappa_2\colon G\lra W_2$.
We will be interested in the induced action of $G$ on $\P(W_1)\times \P(W_2)$. More precisely,
it is the aim of this section to establish a certain property of the instability flag
for a point $(x_1,x_2)\in \P(W_1)\times \P(W_2)$ under the assumption that $x_2\in\P(W_2)$ be
semistable and the chosen polarization be of the form $\O(1,\eta)$ with $\eta\gg 0$. In order to
motivate our result, we need some preparations.
\par
The following discussion is adapted from \cite{Schmitt0}, Section 3.1.
We fix a pair $(B,T)$ in $G$ which consists of a Borel subgroup $B$ and a maximal torus $T\subset B$. 
This defines the following chamber 
$$
\ol{C}:=\bigl\{\,\la\in X_*(T)\,|\, B\subseteq Q_G(\la) \,\bigr\}
$$
in $X_*(T)$. The corresponding rational polyhedral cone in $X_{*,\R}(T)$ will be denoted by
$\ol{C}_\R$. Let $\la_1$,..., $\la_{t^\p}$ be the minimal integral generators of the edges
of $\ol{C}_\R$, so that, in particular,
$$
\ol{C}_\R=\R_{\ge 0}\cdot\la_1+\cdots+\R_{\ge 0}\cdot\la_{t^\p}.
$$
Note that $\ol{C}_\R$ is, in fact, the closure of a Weyl chamber and that every one parameter subgroup of $G$ 
is conjugate to one of $\ol{C}$. Thus, if we
are given a representation $\kappa\colon G\lra \GL(W)$, a point $w\in W\setminus\{0\}$ will be (semi)stable,
if and only if
$$
\mu_\kappa(\la, g\cdot w)(\ge)0,\q \forall g\in G,\ \la\in \ol{C}.
$$
As above, we let
$
{\rm WT}(\kappa,T)
$
be the set of weights of $\kappa$ w.r.t.\ $T$. For a subset $A\subseteq {\rm WT}(\kappa,T)$, we get a
decomposition
$$
\ol{C}_\R=\bigcup_{\chi\in A}C^\chi_A\q\hbox{with}\q
C^\chi_A:=
\bigl\{\,\la\in \ol{C}_\R\,|\,\langle\la,\chi\rangle_\R\ge \langle\la,\chi^\p\rangle_\R\ \forall\chi^\p\in A\bigr\}.
$$
It is easy to see that the intersection $C^\chi_A\cap C^{\chi^\p}_A$, $\chi,\chi^\p\in A$, is a common face of those cones.
If $w\in W\setminus\{0\}$ and $g\in G$, we may take $A:={\rm WT}(g\cdot w,T)$, so that the function 
$$
\mu_\kappa(.,g\cdot w)\colon C^\chi_A\lra \R 
$$
is, on the cone $C^\chi_A$, the \sl linear function \rm $\la\lma \langle\la,\chi\rangle_\R$. 
\par
Let ${\cal S}=\{\, A_1,...,A_s\,\}$ be any subset of the power set of ${\rm WT}(\kappa, T)$.
Since the intersection of two rational polyhedral cones is again a rational polyhedral cone,
it follows easily that we may find a decomposition
$$
\ol{C}_\R=\bigcup_{j\in I({\cal S})} C^{j}_{\cal S}
$$
for an appropriate index set $I({\cal S})$ and rational polyhedral cones $C^{j}_{\cal S}$, $j\in I({\cal S})$,
such that
\begin{itemize}
\item for $j,j^\p\in I({\cal S})$, the intersection $C^j_{\cal S}\cap C^{j^\p}_{\cal S}$ is a common face of these
two cones;
\item for every $A\in {\cal S}$, and every $\chi\in A$, there are indices $j_1,...,j_{s(A,\chi)}\in I({\cal S})$
with $C^\chi_A=\bigcup_{k=1}^{s(A,\chi)} C^{j_k}_{\cal S}$.
\end{itemize}
For our application, we choose
$$
{\cal S}:=\Bigl\{\,{\rm WT}(g\cdot w,T)\,|\, g\in G,w\in W\setminus\{0\}\,\Bigr\}.
$$
Then, for every $g\in G$, every $w\in W\setminus\{0\}$,
and every $j\in I({\cal S})$, the function $\la\lma \mu(\la,g\cdot w)$ is linear on $C^j_{\cal S}$.
We choose additional elements $\la_{t^\p+1},...,\la_{t}\in X_*(T)$, such that
$$
\bigl\{\, \la_1,...,\la_t\,\bigr\}
$$
is the set of elements which occur as minimal integral generators of an edge of a cone $C^j_{\cal S}$, $j\in {\cal S}$.
Thus, a point $w\in W\setminus\{0\}$ will be (semi)stable,
if and only if
$$
\mu_\kappa(\la_i, g\cdot w)(\ge)0,\q \forall g\in G, i\in\{\,1,...,t\,\} .
$$ 
\begin{Prop}
\label{SemStabCrit}
There is a positive rational number $\eta_\infty$, such that for every $\eta> \eta_\infty$,
a point $(x_1,x_2)\in \P(W_1)\times\P(W_2)$ will be (semi)stable w.r.t.\ the linearization
in $\O(1,\eta)$, if and only if {\rm a)} $x_2\in \P(W_2)^{\rm ss}$ and {\rm b)}, for
every one parameter subgroup $\la\colon {\Bbb G}_m(K)\lra G$ with $\mu_{\O_{\P(W_2)}(1)}(\la,x_2)=0$, one has
$$  
\mu_{\O_{\P(W_1)}(1)}(\la,x_1)(\ge)0.
$$
\end{Prop}
\begin{proof}
Set $\widetilde{W}_i:=W^\vee_i$, and let $\widetilde{\kappa}_i\colon G\lra \GL(\widetilde{W}_i)$, $i=1,2$,
be the corresponding representations. This time we look at
$$
{\cal S}:=\Bigl\{\,
{\rm WT}(g\cdot w,T)\,|\, g\in G, w\in\widetilde{W}_1\setminus\{0\}\cup \widetilde{W}_2\setminus\{0\}\,\Bigr\}.
$$
This yields a fan decomposition
$$
\ol{C}_{\R}=\bigcup_{j\in I({\cal S})} C^j_{\cal S},
$$
such that for any $g\in G$, any $w_1\in \widetilde{W}_1\setminus\{0\}$, any $w_2\in \widetilde{W}_2\setminus\{0\}$,
and any $j\in I({\cal S})$, both the function $\la\lma \mu_{\widetilde{\kappa}_1}(.,g\cdot w_1)$
and $\la\lma \mu_{\widetilde{\kappa}_2}(.,g\cdot w_2)$ are linear on $C^j_{\cal S}$. 
Let
$$
\bigl\{\,\la_1,...,\la_t\,\bigr\}\subset X_*(T)
$$
be the set of elements which arise as minimal integral generators of an edge of a cone $C^j_{\cal S}$ for
some $j\in {\cal S}$. 
Set
\begin{eqnarray*}
K_1&:=&\max\Bigl\{\, \mu_{\widetilde{\kappa}_1}(\la_i,g\cdot w_1)\,|\, i=1,...,t,g\in G, w_1\in \widetilde{W}_1
\setminus\{0\}\,\Bigr\}
\\
K_2&:=&\min\Bigl\{\, \mu_{\widetilde{\kappa}_1}(\la_i,g\cdot w_1)\,|\, i=1,...,t,g\in G, w_1\in \widetilde{W}_1
\setminus\{0\}\,\Bigr\}.
\end{eqnarray*}
We first prove that the stated Conditions a) and b) are necessary. Let $x_i=[w_i]$ for 
$w_i\in \widetilde{W}_i\setminus\{0\}$, $i=1,2$, and suppose $\eta> K_1$. If $x_2$ were not semistable, then
there would be an index $i_0\in\{\, 1,...,t\,\}$ and an element $g\in G$ with
$$
\mu_{\O_{\P(W_2)}(1)}(\la_{i_0},x_2)=\mu_{\widetilde{\kappa}_2}(\la_{i_0}, w_2)\le -1,
$$
so that
$$
\mu_{\O(1,\eta)}(\la_{i_0},g\cdot(x_1,x_2))=\mu_{\widetilde{\kappa}_1}(\la_{i_0}, g\cdot w_1)+\eta\cdot
\mu_{\widetilde{\kappa}_2}(\la_{i_0}, g\cdot w_2)\le K_1-\eta<0,
$$
a contradiction. The necessity of Condition b) is now obvious.
\par
Next, suppose $\eta>-K_2$ and that a) and b) are verified. Let $g\in G$ and $i\in\{\,1,...,t\,\}$,
such that $\mu_{\widetilde{\kappa}_2}(\la_i, g\cdot w_2)>0$ (i.e., $\ge 1$).
Then,
$$
\mu_{\widetilde{\kappa}_1}(\la_{i}, g\cdot w_1)+\eta\cdot
\mu_{\widetilde{\kappa}_2}(\la_{i}, g\cdot w_2)>K_2-K_2=0.
$$
If, on the other hand, $\mu_{\widetilde{\kappa}_2}(\la_i, g\cdot w_2)=0$, then Condition b) applies.
Finally, the case $\mu_{\widetilde{\kappa}_2}(\la_i, g\cdot w_2)<0$ is ruled out by Assumption a).
Therefore, we have shown that
$$
\mu_{\O(1,\eta)}(\la_{i},g\cdot (x_1,x_2))(\ge)0\q \forall g\in G,\ i\in\{\,1,...,t\,\},
$$
and that implies that $(x_1,x_2)$ is (semi)stable w.r.t.\ linearization in $\O(1,\eta)$.
The proposition is settled for $\eta_\infty=\max\{\,K_1, -K_2\,\}$.
\end{proof}
This proposition motivates the following result.
\begin{Thm}
\label{InstFlagProduct}
There is a positive rational number $\eta^\p_\infty$, such that for every $\eta> \eta^\p_\infty$
and every point $(x_1,x_2)\in \P(W_1)\times\P(W_2)$ which is not semistable w.r.t.\ the linearization
in $\O(1,\eta)$, but for which $x_2\in \P(W_2)^{\rm ss}$, any instability one parameter
subgroup $\la$ for $(x_1,x_2)$ satisfies
$$  
\mu_{\O_{\P(W_2)}(1)}(\la,x_2)=0.
$$
\end{Thm}
\begin{proof}
In the subsequent considerations, we
will use the same notation and setup as in the proof of Proposition \ref{SemStabCrit}. 
Let $(x_1,x_2)=([w_1],[w_2])$, $w_i\in \widetilde{W}_i\setminus\{0\}$, $i=1,2$, be a point with
$x_2\in \P(W_2)^{\rm ss}$.
Fix elements $g\in G$ and $j\in I({\cal S})$. In the following, we will introduce several constants
which depend both on $g$ and on $j$, but we will mostly omit the letter ``$j$" from the notation.
We do so because $j$ lies a priori in the finite set $I({\cal S})$, and it will be clear that the constants can
be chosen to work for all $j\in I({\cal S})$.
By construction, there exists a $\chi_i\in {\rm WT}(g\cdot w_i, T)$, such that
\begin{equation}
\label{NuFormula}
\nu_{\widetilde{\kappa}_i}(\la, g\cdot w_i)=\frac{\langle\la, \chi_i\rangle_\R}{\| \la\|_G},
\q\hbox{for all $\la\in C^j_{\cal S}$, $i=1,2$}.
\end{equation}
We look at the linear function 
$$
f^g\colon \la\lma \mu_{\widetilde{\kappa}_2}(\la, g\cdot w_2)=\langle\la, \chi_2\rangle_\R
$$
on $C^j_{\cal S}$. By assumption, the function $f^g$ is non-negative on $C^j_{\cal S}$, so that 
$F(g,j):=\{\, f^g=0\,\}\cap C^j_{\cal S}$ is a face of $C^j_{\cal S}$. If $F(g,j)=\varnothing$ and
$\eta>-K_2$, the functions 
$$
\la\lma\mu_{\widetilde{\kappa}_1}(\la,w_1)+\eta\cdot \mu_{\widetilde{\kappa}_2}(\la,w_2) 
$$
and
$$
\la\lma\nu_{\widetilde{\kappa}_1}(\la,w_1)+\eta\cdot \nu_{\widetilde{\kappa}_2}(\la,w_2) 
$$
will be strictly positive on $C^j_{\cal S}$, so that the case will not be interesting for us. The case
$F(g,j)=C^j_{\cal S}$ doesn't have to be considered either (see below).
Thus, we may assume that $F(g,j)$ be a non-empty proper face of $C^j_{\cal S}$.
We choose an affine hyperplane $L\subset X_{*,\R}(T)$ (not containing the origin), such that
$C^j_{\cal S}$ is the cone over the polyhedron
$$
P:=P^j_{\cal S}:=C^j_{\cal S}\cap L.
$$
Then, $P$ is the convex hull of the points, say, $x_1,...,x_m\in L$. After renumbering, we find an index
$1\le i<m$ with 
$$
\{\, x_1,...,x_m\,\}\cap F(g,j)=\{\, x_1,...,x_i\,\}.
$$
Define
$$
P_1:={\rm ConvexHull}(x_1,...,x_i)\q\hbox{and}\q P_2:={\rm ConvexHull}(x_{i+1},...,x_m),
$$
so that $P_1=P\cap F(g,j)$. For any two points $x,y\in X_{*,\R}(T)$, we define
$$
[x,y]:=\bigl\{\,(1-t)\cdot x+t\cdot y\,|\,t\in [0,1]\,\bigr\}.
$$
Clearly,
\begin{equation}
\label{AllOnLine}
P=\bigcup_{x\in P_1, y\in P_2}[x,y].
\end{equation}
For any $\eps>0$, let $P_\eps$ be the set of points in $X_{*,\R}(T)$ which have distance less than $\eps$
to the set $P_1$.
\begin{Claim}
For any $\eps>0$, there exists an $\eta_\eps$, such that for every $\eta>\eta_\eps$
and every $\la\in C^j_{\cal S}$, the
inequality
$$
\nu_{\widetilde{\kappa}_1}(\la,g\cdot w_1)+\eps\cdot \nu_{\widetilde{\kappa}_2}(\la,g\cdot w_2)<0
$$
implies
$$
\la\in P_\eps.
$$
\end{Claim}
We first prove this claim. For this, we introduce the following numbers
\begin{eqnarray*}
\ol{K}^g_1&:=&\min\bigl\{\,\mu_{\widetilde{\kappa}_2}(\la,w_2)\,|\,\la\in P_2\,\bigr\}
\\
\ol{K}^g_2&:=&\max\bigl\{\,\mu_{\widetilde{\kappa}_1}(\la,w_1)\,|\,\la\in P_2\,\bigr\}
\\
\ol{K}^g_3&:=&\min\bigl\{\,\mu_{\widetilde{\kappa}_1}(\la,w_1)\,|\,\la\in P\,\bigr\}
\\
\ol{K}^g_4&:=&\max\bigl\{\,\| x-y\|_G\,|\,x\in P_1, y\in P_2\,\bigr\}.
\end{eqnarray*}
We also define
\begin{equation}
\label{Eps101}
\eps^\p:=\frac{\eps}{\ol{K}^g_4}.
\end{equation}
Note that $\ol{K}^g_1$ is a strictly positive number.
Suppose we are given $x\in P_1$ and $y\in P_2$. We define the linear function
\begin{eqnarray*}
h^\eta_{x,y}\colon [0,1] &\lra &\R
\\
t &\lma& \mu_{\widetilde{\kappa}_1}\bigl((1-t)x+ty, g\cdot w_1\bigr)+\eta\cdot 
\mu_{\widetilde{\kappa}_2}\bigl((1-t)x+ty, g\cdot w_2\bigr).
\end{eqnarray*}
The slope of this function is
\begin{eqnarray}
h^\eta_{x,y}(1)-h^\eta_{x,y}(0)&\stackrel{\mu_{\widetilde{\kappa}_1}(x,g\cdot w_2)=0}{=}&
\mu_{\widetilde{\kappa}_1}\bigl(y, g\cdot w_1\bigr)+\eta\cdot \mu_{\widetilde{\kappa}_2}\bigl(y, g\cdot w_2\bigr)
-\mu_{\widetilde{\kappa}_1}\bigl(x, g\cdot w_1\bigr)
\\ 
\label{TestEQ1}
&\ge & \ol{K}^g_3+\eta\cdot \ol{K}^g_1-\ol{K}^g_2.
\end{eqnarray}
Moreover, we compute
\begin{eqnarray}
h^\eta_{x,y}(\eps^\p)&=&
(1-\eps^\p)\cdot \mu_{\widetilde{\kappa}_1}\bigl(x, g\cdot w_1\bigr)+\eps^\p\cdot\mu_{\widetilde{\kappa}_1}\bigl(y, g\cdot w_1\bigr)+\eta\cdot\eps^\p\cdot\mu_{\widetilde{\kappa}_2}\bigl(y, g\cdot w_2\bigr)
\\ 
\label{TestEQ2}
&\ge & \ol{K}^g_3+\eta\cdot\eps^\p\cdot \ol{K}^g_1.
\end{eqnarray}
Choose
$$
\eta>\eta^g_{\eps^\p}:=
\max\left\{\,\frac{\ol{K}^g_2-\ol{K}^g_3}{\ol{K}^g_1}, \frac{-\ol{K}^g_3}{\eps^\p\cdot \ol{K}^g_1}\,\right\}.
$$
Suppose
$$
\nu_{\widetilde{\kappa}_1}(\la,g\cdot w_1)+\eps\cdot \nu_{\widetilde{\kappa}_2}(\la,g\cdot w_2)<0.
$$
By (\ref{AllOnLine}), we find points $x\in P_1$, $y\in P_2$, and $t_0\in [0,1]$
with
$$
\la=(1-t_0)x+t_0y.
$$
Since
$$
\nu_{\widetilde{\kappa}_1}(\la,g\cdot w_1)+\eps\cdot \nu_{\widetilde{\kappa}_2}(\la,g\cdot w_2)
=h^\eta_{x,y}(t_0),
$$
(\ref{TestEQ1}) and (\ref{TestEQ2}) imply $t_0<\eps^\p$.
One computes
$$
\| (1-t_0)x+t_0y-x\|_G=t_0\cdot\| y-x\|_G<\eps^\p\cdot \| y-x\|_G
\le \eps^\p\cdot \ol{K}^g_4\stackrel{(\ref{Eps101})}{=}\eps.
$$
Thus, $\la\in P_\eps$, and the claim is settled.
\par
Next, we define, for $i=1,2$, the differentiable function
\begin{eqnarray*}
{\frak N}^g_i\colon \R\times P_1\times P_2 &\lra &\R
\\
(t,x,y) &\lma& \frac{\langle\,(1-t)\cdot x+t\cdot y,\chi_i \,\rangle_\R}{\|(1-t)\cdot x+t\cdot y\|_G}
\end{eqnarray*}
(which coincides with $\nu_{\widetilde{\kappa}_i}\bigl((1-t)\cdot x+t\cdot y, g\cdot w_i\bigr)$
for $0\le t\le 1$, $i=1,2$, by (\ref{NuFormula})) as well as
\begin{eqnarray*}
{\frak D}^g_i\colon \R\times P_1\times P_2 &\lra &\R
\\
(t,x,y) &\lma& \Bigl(\frac{\partial}{\partial t}{\frak N}_i^g\Bigr)(t,x,y).
\end{eqnarray*}
\begin{Claim}
The function ${\frak D}^g_2$ is strictly positive on $\{0\}\times P_1\times P_2$.
\end{Claim}
We now prove this claim. Given $x$ and $y$, we have the function 
\begin{eqnarray*}
N^g_{2;x,y}\colon \R &\lra &\R
\\
t &\lma& {\frak N}^g_i(t,x,y) = \frac{\mu(t)}{n(t)}
\end{eqnarray*}
with
$$
\mu(t):= 
\langle\,(1-t)\cdot x+t\cdot y,\chi_i \,\rangle_\R
\stackrel{(\ref{NuFormula})}{=} t\cdot \mu_{\widetilde{\kappa}_2}\bigl(y, g\cdot w_2\bigr)
\ \hbox{and}\ n(t):=\| (1-t)x+ty\|_G.
$$
We see
$$
{\frak D}^g_2(0,x,y)=(N^g_{2;x,y})^\p(0)=\frac{\mu^\p(0)\cdot n(0)-\mu(0)\cdot n^\p(0)}{n(0)^2}=
\frac{\mu_{\widetilde{\kappa}_2}\bigl(y, g\cdot w_2\bigr)}{n(0)}>0.
$$
Thus, the claim is established.
\par
Because of compactness, there is a number $\eps_0^\p(g)>0$, such that ${\frak D}^g_2$ is strictly
positive on $[0,\eps_0^\p(g)]\times P_1\times P_2$. We define the additional constants
\begin{eqnarray*}
\ol{K}^g_5&:=&\min\bigl\{\,{\frak D}^g_2(t,x,y)\,|\,(t,x,y)\in [0,\eps_0^\p(g)]\times P_1\times P_2\,\bigr\}
\\
\ol{K}^g_6&:=&\min\bigl\{\,{\frak D}^g_1(t,x,y)\,|\,(t,x,y)\in [0,\eps_0^\p(g)]\times P_1\times P_2\,\bigr\}.
\end{eqnarray*}
Choose an $\eps(g)$, such that $P_{\eps(g)}$ is contained in the image of the map 
\begin{eqnarray*}
[0,\eps_0^\p(g)]\times P_1\times P_2 &\lma & P
\\
(t,x,y) &\lma& (1-t)\cdot x+ t\cdot y.
\end{eqnarray*}
\begin{Claim}
Suppose 
$$
\eta> \eta^\p(g):=-\frac{\ol{K}^g_6}{\ol{K}^g_5},
$$
then
$$
\min\Bigl\{\,\nu_{\widetilde{\kappa}_1}(\la,g\cdot w_1)+\eta\cdot \nu_{\widetilde{\kappa}_2}(\la,g\cdot w_2)\,\Big|\,
\la\in P_{\eps(g)}\,\Bigr\}
$$
is taken on at a point $\la\in P_1$.
\end{Claim}
To prove the claim, we first note that, by (\ref{AllOnLine}) and our choice of $\eps(g)$, we may write 
$$
\la=(1-t_0)x+t_0y\q\hbox{for some $t_0\le \eps^\p_0(g)$}.
$$
We have to look at the function
$$
N^g_{x,y}(t):={\frak N}_1^g(t,x,y)+\eta\cdot {\frak N}_2^g(t,x,y).
$$
By assumption
$$
(N^g_{x,y})^\p(t)= {\frak D}_1^g(t,x,y)+\eta\cdot {\frak D}_2^g(t,x,y)>0\q\hbox{for $0\le t\le \eps^\p(g)$}.
$$
Therefore, $N^g_{x,y}$ is strictly monotonously growing on $[0,\eps^\p(g)]$, whence
$$
\nu_{\widetilde{\kappa}_1}(x,g\cdot w_1)+\eta\cdot \nu_{\widetilde{\kappa}_2}(x,g\cdot w_2)
=N^g_{x,y}(0)\le N^g_{x,y}(t_0)
=\nu_{\widetilde{\kappa}_1}(\la,g\cdot w_1)+\eta\cdot \nu_{\widetilde{\kappa}_2}(\la,g\cdot w_2)
$$
with equality if and only if $t_0=0$, i.e., $\la=x\in P_1$. This clearly settles the claim.
\par
We are finally in position to prove the theorem. Since the translates of $\ol{C}_{\R}$ under the Weyl group
cover $X_{*,\R}(T)$, it is clear from the sketch of the proof of Theorem \ref{kempfi} that
\begin{eqnarray*}
\min_{\gamma\in \Gamma}\min\left\{\, \frac{l_\eta^\gamma(\la)}{\|\la\|_G}\,\Big|\, \la\in X^*_\R(T)\,\right\}
&=&
\min_{\gamma\in \Gamma}\min\left\{\, \frac{l_\eta^\gamma(\la)}{\|\la\|_G}\,\Big|\, \la\in \ol{C}_\R\,\right\}
\\
&=&
\min_{\gamma\in \Gamma}\min_{j\in I({\cal S})}
\min\left\{\,\frac{l_\eta^\gamma(\la)}{\|\la\|_G}\,\Big|\, \la\in C^j_{\cal S}\,\right\}.
\end{eqnarray*}
Here, $l_\eta^g$ is the function formed w.r.t.\ a representation $S^{r_1}(\widetilde{\kappa}_1)\otimes
S^{r_2}(\widetilde{\kappa}_2)$, such that $r_1,r_2\in \Z_{>0}$ satisfy $r_1/r_2=1/\eta$.
Now, choose
$$
\eps\le \min_{\gamma\in \Gamma}\{\, \eps(\gamma)\,\}
$$
and 
$$
\eta^\p_\infty>\max_{\gamma\in \Gamma}\{\,\eta^\gamma_\eps, \eta^\p(\gamma)\,\}.
$$
Suppose $\eta>\eta^\p_\infty$, and let $\la$ be an instability one parameter subgroup for $(x_1,x_2)$. Then, we find
a $g\in G$, a $\gamma_0\in G$, and a $j_0\in I({\cal S})$, such that
\begin{equation}
\label{NearTheEnd}
\hbox{$\la_0:=g\cdot \la\cdot g^{-1}\in C^{j_0}_{\cal S}$, and ${\rm WT}(\gamma_0\cdot w_i, T)=
{\rm WT}(g\cdot w_i, T)$},\ i=1,2.
\end{equation}
Being an instability one parameter subgroup implies
\begin{eqnarray*}
\nu_{\widetilde{\kappa}_1}(\la_0,\gamma_0\cdot w_1)+ \eta
\nu_{\widetilde{\kappa}_2}(\la_0,\gamma_0\cdot w_2)
&=&\min\bigl\{\, \nu_{\widetilde{\kappa}_1}(\la,\gamma_0\cdot w_1)+ \eta
\nu_{\widetilde{\kappa}_2}(\la,\gamma_0\cdot w_2)\,|\, \la\in C^{j_0}_{\cal S}\,\bigr\}
\\
&=&\min\bigl\{\, \nu_{\widetilde{\kappa}_1}(\la,\gamma_0\cdot w_1)+ \eta
\nu_{\widetilde{\kappa}_2}(\la,\gamma_0\cdot w_2)\,|\, \la\in P^{j_0}_{\cal S}\,\bigr\}
\\
&<&0.
\end{eqnarray*}
As remarked at the beginning $F(\gamma_0, j_0)$ must be a non-empty face of $C^{j_0}_{\cal S}$.
If it is all of $C^{j_0}_{\cal S}$, then 
$$
0=\mu_{\widetilde{\kappa}_2}(\la_0,\gamma_0\cdot w_2)\stackrel{(\ref{NearTheEnd})}{=}
\mu_{\widetilde{\kappa}_2}(\la_0,g\cdot w_2)=\mu_{\widetilde{\kappa}_2}(\la, w_2)
$$
which is what we want. Otherwise, by the first claim we made and our choice of $\eta$, we
have $\widetilde{\la}_0\in P_\eps$. Here, $\widetilde{\la}_0$ is the point of intersection
of the ray $\R_{\ge 0}\cdot \la_0$ with $P^\chi_A$. But then, the third claim implies
$\widetilde{\la}_0\in P_1$, i.e., 
$$
\nu_{\widetilde{\kappa}_2}(\la_0,\gamma_0\cdot w_2)=
\nu_{\widetilde{\kappa}_2}(\widetilde{\la}_0,\gamma_0\cdot w_2)=0, 
$$
whence again
$\mu_{\widetilde{\kappa}_2}(\la_0,\gamma_0\cdot w_2)=0$.
The proof of the theorem is now complete.
\end{proof}
\paragraph{Simplification of the Semistability Concept for Decorated Sheaves. ---}
Here, we work with the group $G=\GL_r(\C)$. For given non-negative integers $a$, $b$, and $c$,
and the complex vector space $V:=\C^r$,
we have the representation $\kappa_{a,b,c}$ of $G$ on 
$W:=(V^{\otimes a})^{\oplus b}\otimes (\bigwedge^rV)^{\otimes -c}$. Again, $\widetilde{\kappa}$
will be the corresponding contragredient representation of $\GL_r(\C)$ on
$\widetilde{W}:=W^\vee$. We carry out the discussion from the beginning of the Section ``The Instability Flag
in a Product" for the representation $\widetilde{\kappa}$. This provides us with one parameter subgroups
$\la_1,...,\la_t$. Furthermore, the one parameter subgroup $\la_i$ yields the weighted flag
$$
\Bigl(\, V_i^\bullet\colon\q 0\subsetneq V_1^i\subsetneq\cdots\subsetneq V^i_{s_i}\subsetneq \C^r,
\ul{\alpha}_i=(\alpha^i_1,...,\alpha^i_{s_i})\,\Bigr)
$$
in $V$, $i=1,...,t$.  We define the finite set
$$
{\cal T}:=\Bigl\{\,(\dim V^i_1,...,\dim V^i_{s_i}; \alpha^i_1,...,\alpha^i_{s_i})\,\Big|\,i=1,...,t\,\Bigr\}.
$$
Then, one has the following result.
\begin{Prop}
\label{ConeDecompII}
Let $(\E,\phi)$ be a torsion free sheaf of rank $r$ with a decoration of type $(a,b,c;\L)$,
and $\delta\in\Q[x]$ a positive polynomial of degree at most $\dim X-1$. Then,
$(\E,\phi)$ is $\delta$-(semi)stable, if and only one has
$$
M(\E,\ul{\alpha})+\delta\cdot \mu(\E^\bullet,\ul{\alpha};\phi)(\succeq)0
$$
for every weighted filtration $(\E,\ul{\alpha})$ with
$$
(\rk\E_1,...,\rk\E_s;\alpha_1,...,\alpha_s)\in {\cal T}.
$$
\end{Prop}
\begin{proof}
This follows immediately from the general formalism. More details are contained in \cite{Schmitt0},
especially Section 3.1.
\end{proof}
Examples of how this proposition simplifies the semistability concept may be found in Chapter 3 of \cite{Schmitt0}.
In fact, it is a very important finiteness result which we will use below.
\subsection{Principal Bundles}
Let $U$ be a smooth algebraic variety and $G$ a reductive algebraic group over the field of complex numbers. 
Suppose we are given a principal
$G$-bundle ${\cal P}$ over $U$. If $F$ is an algebraic variety and $\alpha\colon G\times F\lra F$ is an 
action of $G$ on $F$, then we may form the geometric quotient
$$
{\cal P}(F,\alpha):= \bigl({\cal P}\times F\bigr)/G
$$ 
w.r.t.\ the action $(p,f)\cdot g:= (p\cdot g, g^{-1}\cdot f)$ for all $p\in {\cal P}$, $f\in F$, and $g\in G$. Note that
${\cal P}(F,\alpha)$ is a fibre space with fibre $F$ over $U$ which is
locally trivial in the \'etale topology. An important 
special case arises when we look at the action $\frak c\colon G\times G\lra G$, $(g,h)\lma g\cdot h\cdot g^{-1}$,
of $G$ on itself by conjugation. Then, the associated fibre space ${\cal G}({\cal P}):={\cal P}(G,\frak c)\lra U$ is 
a reductive group scheme over $U$, and, for any pair $(F,\alpha)$ as above, we obtain an induced action
$$
a\colon\ {\cal G}({\cal P})\times_U {\cal P}(F,\alpha)\lra {\cal P}(F,\alpha).
$$
If $W$ is a vector space and $\kappa\colon G\lra \GL(W)$ is a representation
and $\ol{\alpha}\colon G\times \P(W)\lra \P(W)$ is the corresponding action, 
we set ${\cal P}_\kappa:={\cal P}(W,\kappa)$ and $\ol{\cal P}_\kappa:={\cal P}(\P(W),\ol{\alpha})$.
Note that the formation of ${\cal P}(F,\alpha)$ commutes with base change. For additional information, we refer
the reader to \cite{Serre}.
\begin{Ex}
Suppose $U\subset X$ is a big open subset and that $E$ is a vector bundle on $U$ with
typical fibre $V$. Let $\la\colon \C^*\lra \SL(V)$ be a non-trivial one parameter subgroup and
$(V^\bullet,\ul{\alpha})$ the weighted flag of $\la$. Then, by construction,
$Q_{\SL(V)}(\la)$ and $Q:=Q_{\GL(V)}(\la)$ are the $\SL(V)$- and the $\GL(V)$-stabilizer of $V^\bullet$, respectively. 
Finally, 
${\cal P}:=\ul{\rm Isom}(V\otimes\O_U,E)$ 
is the principal $\GL(V)$-bundle associated with $E$. We have the natural left
action $\ol{\beta}$ of $\GL(V)$ on the projective variety $\GL(V)/Q$ and
$$
{\cal Q}:={\cal P}/Q\cong {\cal P}(\GL(V)/Q,\ol{\beta}).
$$
Now, suppose $U^\p\subseteq U$ is a big open subset. Then, giving a section
$\sigma_{U^\p}\colon U^\p\lra  {\cal Q}$ is the same as giving a filtration
$0\subsetneq E_1\subsetneq\cdots\subsetneq E_s\subsetneq E_{|U^\p}$ of $E_{|U^\p}$
by subbundles with $\rk E_i=\dim V_i$, $i=1,...,s$.
As we have seen in Example \ref{CharComp}, the one parameter subgroup $\la$ gives rise
to a character $\chi_\la$ of $Q$. Thus, we obtain the line bundle 
$$
\L_{\chi_\la}:={\cal Q}(\C, \chi).
$$
We claim
\begin{equation}
\label{DegComp}
\deg\L_{\chi_\la}=\sum_{i=1}^s \alpha_i\bigl(\deg E\cdot\rk E_i - \deg E_i\cdot\rk E\bigr).
\end{equation}
For this, let $\ul{v}=(v_1,...,v_r)$ be a basis consisting of eigenvectors for $\la$,
i.e., $\la(z)\cdot v_i=z^{\gamma_i}\cdot v_i$, $i=1,...,r$, and we assume
$\gamma_1\le\gamma_2\le\cdots\le \gamma_r$. Denote by $T$ the maximal torus defined by $\ul{v}$. First, consider the case that $s=1$
and $\alpha_1=1$, then
$$
\ul{\gamma}:=(\gamma_1,...,\gamma_r)=\gamma^{(j)}:=\bigl(\underbrace{j-r,...,j-r}_{j\times},
\underbrace{j,...,j}_{(r-j)\times} \bigr),\q j:=\dim V_1=\rk E_1.
$$
Now, Example \ref{CharComp} implies 
$$
\L_{\chi_\la}\cong \det(E_1)^{\otimes (\rk E_1-r)}\otimes
\det(E/E_1)^{\otimes \rk E_1},
$$
and one easily checks Equation (\ref{DegComp}) in this case.
Next, assume that $\alpha_1,...,\alpha_s$ be integers.
Letting $\la^{(j)}$ be the one parameter subgroup defined by the weight vector $\gamma^{(j)}$
and the basis $\ul{v}$, $j=1,...,r-1$,  
we see that 
$$
\la=\sum_{i=1}^s \alpha_i\cdot \la^{(\dim V_i)}\q\hbox{in $X_*(T)$},
$$
whence $\chi_\la=\sum_{i=1}^s \alpha_i\cdot \chi_{\la^{(\dim V_i)}}$ in $X^*(T)$ and
$$
\L_{\chi_\la}\cong \bigotimes_{i=1}^s \det(\L_{\la^{(\dim V_i)}})^{\otimes \alpha_i}.
$$
Thus, Equation (\ref{DegComp}) follows easily. In the remaining case, we choose a positive
integer $q$, such that $q\alpha_i\in \Z_{>0}$, $i=1,...,s$, and use the case
we have just treated together with the fact $\L_{\chi_{q\la}}\cong \L_{\chi_\la}^{\otimes q}$.
\qed
\end{Ex}
\paragraph{Parabolic Subgroup Schemes. ---}
Let $S$ be any scheme and
suppose ${\cal G}_S\lra S$ is a reductive group scheme over $S$. A subgroup ${\cal Q}_S\subset {\cal G}_S$ is called
a \it parabolic subgroup\rm, if it is smooth over $S$ and, for any geometric point $s$ of $S$, the quotient
${\cal G}_{S,s}/{\cal Q}_{S,s}$ is proper. The functor
\begin{eqnarray*}
\ul{\rm Par}({\cal G}_S)\colon\ \ul{\rm Schemes}_S &\lra& \ul{\rm Sets}
\\
(T\lra S) &\lma& \Bigl\{\, \hbox{Parabolic subgroups of ${\cal G}_S\times_S T$}\,\Bigr\}
\end{eqnarray*}
is then representable by an $S$-scheme ${\cal P}ar({\cal G}_S)$. For the details, we refer the reader to \cite{Dem}.
\begin{Ex}
\label{ParRep}
Let $G$ be a complex reductive group and $\cal P$ be a principal $G$-bundle over the variety $U$. 
Denote by $\frak P$ the
set of conjugacy classes of parabolic subgroups in $G$ and pick for each class $\frak p\in \frak P$ a representative
$Q_\frak p$. Then,
$$
{\cal P}ar\bigl({\cal G}({\cal P})\bigr)\cong \bigsqcup_{\frak p\in\frak P} {\cal P}/Q_\frak p.
$$
\end{Ex}
\paragraph{Sections in Associated Projective Bundles. ---}
Let ${\cal P}\lra U$ be a principal bundle as before. Suppose we are given a representation $\kappa\colon
G\lra \GL(W)$. This yields an action $\alpha\colon G\times \P(W^\vee)\lra \P(W^\vee)$ and a linearization
$\ol{\alpha}\colon G\times \O_{\P(W^\vee)}(1)\lra \O_{\P(W^\vee)}(1)$ of this action. 
Let $\P^{\rm ss}\subset \P(W^\vee)$ be the open subset of semistable points.
\begin{Prop}[Ramanan-Ramanathan] 
\label{RamRamII}
Assume that $U$ is a big open subset of the manifold $X$. Let  
$\sigma\colon U\lra {\cal P}(\P(W^\vee),\alpha)$ be a section, and $\L_\sigma$ the pullback --- via
$\sigma$ --- of
the line bundle ${\cal P}(\O_{\P(W^\vee)}(1),\ol{\alpha})\lra {\cal P}(\P(W^\vee),\alpha)$ to $U$. 
If $\sigma(\eta)\in
{\cal P}(\P^{\rm ss},\alpha)$, then
$$
\deg(\L_\sigma)\q\ge\q 0.
$$
Here, $\eta$ is the generic point of $U$.
\end{Prop}
\begin{proof} This is Proposition 3.10, i), in \cite{RamRam}.
\end{proof}
Occasionally, we say that the section $\sigma$ is \it generically semistable\rm, if $\sigma(\eta)\in
{\cal P}(\P^{\rm ss},\alpha)$.
\section{The Proof of the Main Theorem}
As in the introduction, we fix a Hilbert polynomial $P$ and a line bundle $\L$ on $X$ of degree,
say, $d$.
We recall some notation. Let $a$, $b$, and $c$ be non-negative integers, $r$ the rank of
sheaves on $X$ with Hilbert polynomial $P$, and fix a complex vector space $V$ of dimension
$r$. Define
$$
V_{a,b,c}:= \bigl(V^{\otimes a}\bigr)^{\oplus b}\otimes \bigl(\bigwedge^r V\bigr)^{\otimes -c},
$$ 
and let $\kappa:=\kappa_{a,b,c}\colon \GL(V)\lra \GL(V_{a,b,c})$ be the corresponding representation,
and $\alpha\colon \GL(V)\times Y\lra Y$, $Y:=\P(V_{a,b,c})$ the induced action. 
Finally, denote by ${\Bbb P}^{\rm ss}_{a,b,c}$
the set of $\SL(V)$-semistable points in $\P(V_{a,b,c})$.
\subsection{Step 1}
\begin{Thm}
\label{SemStabRed}
There is a constant $\delta_0$, depending only on the input data $a$, $b$, $c$, $d$, and $P$,
such that for every stability parameter $\delta$ of degree exactly $\dim(X)-1$ the leading
coefficient of which is larger than $\delta_0$ and every 
$\delta$-semistable torsion free sheaf $(\E,\phi)$ with a decoration of type $(a,b,c;\L)$
and $P(\E)=P$, the following holds true: Denote by $\eta$ the generic point of $X$ and by $K$ its residue field, 
and choose a trivialization
$\E\otimes_{\O_X} \O_{X,\eta}\cong V\times_{\Spec(\C)} \Spec(K)$.
Then, the point $\widetilde{\sigma}_\eta\in \P(V_{a,b,c})\times_{\Spec(\C)} \Spec(K)$ 
defined by $\phi$ and the chosen trivialization lies
in $\P^{\rm ss}_{a,b,c}\times_{\Spec(\C)} \Spec(K)$. 
\end{Thm}
\begin{proof}
Define ${\cal P}:=\ul{\rm Isom}(V\otimes\O_U,\E_{|U})$ as the principal $\GL(V)$-bundle
associated with $\E_{|U}$ over the maximal open subset $U$ over which $\E$ is locally free. 
The group $\GL(V)$ acts on $\SL(V)$ and $\GL(V)$ by conjugation $\frak c$, and we let
${\cal SL}({\cal P})\subset {\cal GL}({\cal P})$ be the corresponding group schemes over $U$.
We fix an algebraic
closure ${\Bbb K}$ of $K$. A trivialization as chosen in the statement of the theorem
is equivalent to a trivialization ${\cal P}\times_U \Spec(K)\cong \GL(V)\times_{\Spec(\C)} \Spec(K)$.
The latter identification will induce trivializations of the objects introduced below.
Define
$$
\ol{Y}:= {\cal P}(Y,\alpha)\times_U \Spec({\Bbb K})\cong Y\times_{\Spec(\C)} \Spec({\Bbb K}),
$$
$G:=\SL(V)$, and
$$
\ol{G}:= {\cal SL}({\cal P})\times_U \Spec({\Bbb K})\cong G\times_{\Spec(\C)} \Spec({\Bbb K}).
$$
Finally, set $L:=\O_Y(1)$, and $\ol{L}:=\O_{\ol{Y}}(1)\cong L\times_{\Spec(\C)} \Spec({\Bbb K})$.
\par
Next, we remind the reader of Proposition 1.14 in Chapter 1.4 of \cite{GIT}:
\begin{Prop}[Mumford]
\label{Mummy}
The set of $\ol{G}$-semistable points in $\ol{Y}$ w.r.t.\ the linearization in $\ol{L}$ is given as
$$
Y^{\rm ss}(L)\times_{\Spec(\C)} \Spec({\Bbb K}).
$$
Here, $Y^{\rm ss}(L)$ is the set of $G$-semistable points in $Y$ w.r.t.\ linearization in $L$.
\end{Prop}
Now, let ${\cal NC}\subset V^\vee_{\ul{a},\ul{b},\ul{c}}$ be the cone of $\SL(V)$-unstable points.
Recall that, over the open subset $\widetilde{U}$ over which $\phi$ is surjective, 
we are given a section $\widetilde{\sigma}_{\widetilde{U}}\colon \widetilde{U}\lra 
\ol{\cal P}_\kappa$. Set
$\widetilde{\sigma}_\eta:=\widetilde{\sigma}_{\widetilde{U}}\times_{\widetilde{U}}\Spec({\Bbb K})\in \ol{Y}$. 
The negation of the assertion of the theorem is,
by Proposition~\ref{Mummy},
$$
\widetilde{\sigma}_\eta\in\bigl({\cal NC}/\C^*\bigr)\times_{\Spec(\C)} \Spec({\Bbb K}).
$$
Our first step toward the proof will be an application of Kempf's rationality theorem \ref{rationality}.
For this, let $T\subset\GL(V)$ be a maximal torus. We may choose a basis of $V$, such that $T$ becomes
the subgroup of diagonal matrices. Then, we define the pairing $(.,.)^*\colon X^*_{\Bbb R}(T)\times X^*_{\Bbb R}(T)
\lra \R$ as in Section \ref{GIT} (The Instability Flag). Now, $\ol{T}:= T\times_{\Spec(\C)} \Spec({\Bbb K})$
is a maximal torus in $\GL(V)\times_{\Spec(\C)} \Spec({\Bbb K})$ with $X^*(\ol{T})=X^*(T)$, 
and its intersection $T_{\ol{G}}$ with
$\ol{G}$ is a maximal torus in that group. The induced pairing $(.,.)^*_G$ on $X_\R^*(T_{\ol{G}})$ fulfills
the requirements of Theorem~\ref{rationality}. If we assume that $\widetilde{\sigma}_\eta$ be unstable, then
there is an instability one parameter subgroup $\la\colon {\Bbb G}_m({\Bbb K})\lra \ol{G}$ which
defines a weighted flag $(0\subsetneq \ol{V}_1\subsetneq\cdots\subsetneq \ol{V}_s\subsetneq \ol{V},\ul{\alpha})$
in $\ol{V}:=V\otimes_\C {\Bbb K}$. The resulting parabolic subgroup $Q_{\ol{G}}(\widetilde{\sigma}_\eta)$ is 
defined over $K$, i.e., it comes from a parabolic subgroup $Q_{K}(\widetilde{\sigma}_\eta)$ of ${\cal SL}({\cal P})\times_U
\Spec(K)$. The parabolic subgroup $Q_{K}(\widetilde{\sigma}_\eta)$, in turn, corresponds to a point
$$
\Spec(K)\lra {\cal P}ar({\cal SL}({\cal P}))={\cal P}ar({\cal GL}({\cal P}))\stackrel{{\rm \ref{ParRep}}}{\cong}
{\cal P}/Q_{\frak p},
$$
for the appropriate conjugacy class ${\frak p}\in {\frak P}$ of parabolic subgroups of $\GL(V)$.
This point may be extended to a section
$$
U^\p\lra 
{\cal P}/Q_{\frak p}
$$
over a non-empty open subset $U^\p\subseteq U$. In fact, we may assume $U^\p$ to be big. This is because
$X$ and ${\cal P}/Q_{\frak p}$ are smooth projective varieties, so that any rational map
$X\dasharrow {\cal P}/Q_{\frak p}$ extends to a big open subset. The group $Q_{\frak p}$ is the stabilizer of a unique flag $V^\bullet\colon 0\subsetneq V_1\subsetneq\cdots \subsetneq V_s\subsetneq V$. Therefore, ${\cal P}/Q_{\frak p}$ identifies with the bundle of flags 
in the fibres of $\E_{|U}$ having the same dimension vector as the flag $V^\bullet$, so that,
from $U^\p\lra 
{\cal P}/Q_{\frak p}$, we get a filtration
$$
0\subsetneq {\cal E}^\p_1\subsetneq \cdots \subsetneq {\cal E}^\p_s\subsetneq \E_{|U^\p}
$$
of $\E_{|U^\p}$ by subbundles. We define $\E_i$ as the saturation of $j_*(\E^\p_i)\cap {\cal E}$
in $\E$, $i=1,...,s$, $j\colon U^\p\lra X$ being the inclusion, i.e.,
we get a filtration 
$$
\E^\bullet\colon 0\subsetneq \E_1\subsetneq \cdots \subsetneq \E_s\subsetneq \E
$$
of $\E$ by saturated subsheaves. 
For the weighted filtration $(\E^\bullet, \ul{\alpha})$, we clearly find
$$
\mu\bigl(\E^\bullet, \ul{\alpha}; \phi\bigr) < 0.
$$
Thus, $\delta\cdot \mu(\E^\bullet, \ul{\alpha}; \phi)$ is a negative
polynomial of degree exactly $\dim(X)-1$. 
\par
Note that we obtain, in fact, an even stronger rationality theorem. As remarked above, the group $Q_{\frak p}$ is the stabilizer
of a unique flag $V^\bullet\colon 0\subsetneq V_1\subsetneq\cdots \subsetneq V_s\subsetneq V$, and
the weighted filtration $(\E^\bullet, \ul{\alpha})$ defines a reduction of the structure group of ${\cal P}$ to $Q_{\frak p}$. If we start our arguments with a trivialization of the induced 
$Q_{\frak p}$-bundle ${\cal Q}$
over the generic point $\eta$, then we get $\ol{V}^\bullet$ with $\ol{V}_i=V_i\otimes_\C {\Bbb K}$, $i=1,...,s$,
as the instability flag. One may use the weighted flag 
$(V^\bullet, \ul{\alpha})$ to define
a one parameter subgroup $\la\colon \C^*\lra G$ (which, indeed, is an instability subgroup). 
Then, $\la$ defines also a flag 
$W^\bullet\colon\ 0\subsetneq W_1\subsetneq \cdots\subsetneq W_{t}
\subsetneq W$ in $W:=V^\vee_{\ul{a},\ul{b},\ul{c}}$, 
and the parabolic subgroup $Q_{G}(\la)\subset Q_{\GL(V)}(\la)=Q_{\frak p}$ fixes this flag. 
Recall that we are given a reduction of the structure group of ${\cal P}_{|U^\p}$ 
to $Q_{\GL(V)}(\la)=Q_{\frak p}$. Therefore, the flag $W^\bullet$ gives rise to a filtration 
$$
0\subsetneq {\cal F}^\vee_1\subsetneq \cdots\subsetneq {\cal F}^\vee_{t}\subsetneq \E_{|U^\p;\ul{a},\ul{b},\ul{c}}^\vee
$$
by subbundles.
Define 
$$
j_0:=\min\{\, j=1,...,t+1\,|\, \P({\cal F}_j)\hbox{ contains the image of $\widetilde{\sigma}_{U}$}\,\}.
$$
Let $\L^\p\subset \L$ be the image of $\phi$. Then, over a big open subset $U^{\p\p}\subseteq U$,
we have $\L^\p_{|U^{\p\p}}\cong \L_{|U^{\p\p}}(-D)$ for an effective divisor $D$. Thus,
$\phi_{|U^{\p\p}}\colon \E_{a,b,c|U^{\p\p}}\lra \L_{|U^{\p\p}}(-D)$ defines a morphism
$$
\ol{\sigma}_{U^{\p\p}}\colon U^{\p\p}\lra \P(\E_{a,b,c|U^{\p\p}})
$$
with
$$
{\ol{\sigma}_{U^{\p\p}}}^*\bigl(\O_{\P(\E_{a,b,c|U^{\p\p}})}(1)\bigr)\cong\L_{|U^{\p\p}}(-D).
$$
By our choice of $j_0$, 
$\ol{\sigma}_{U^{\p\p}}$ factorizes over $\P({\cal F}_{j_0|U^{\p\p}})$,
and, again,
$$
{\ol{\sigma}_{U^{\p\p}}}^*\bigl(\O_{\P({\cal F}_{j_0|U^{\p\p}})}(1)\bigr)\cong\L_{|U^{\p\p}}(-D).
$$
Now, the  surjective linear map $W_{j_0}\lra W_{j_0}/W_{j_0-1}$ is, in fact, a morphism
of $Q_{\frak p}$-modules.
Over a big open subset $U^{\p\p\p}\subseteq U^{\p\p}$, the image of 
$$
\bigl({\cal F}^\vee_{j_0|U^{\p\p\p}}/{\cal F}^\vee_{j_0-1|U^{\p\p\p}}\bigr)^\vee \subset {\cal F}_{j_0|U^{\p\p\p}}
\lra \L_{|U^{\p\p\p}}(-D)
$$
is of the form $\L_{|U^{\p\p\p}}(-(D_{|U^{\p\p\p}}+D^\p))$ for some effective divisor $D^\p$. Therefore, we get
a morphism
$$
\ol{\sigma}^{\p\p}\colon U^{\p\p\p}\lra 
\P\bigl(({\cal F}^\vee_{j_0|U^{\p\p\p}}/{\cal F}^\vee_{j_0-1|U^{\p\p\p}})^\vee\bigr)
$$
with
$$
{\ol{\sigma}^{\p\p}}^*\bigl(\O_{\P(({\cal F}^\vee_{j_0|U^{\p\p\p}}/{\cal F}^\vee_{j_0-1|U^{\p\p\p}})^\vee)}(1)\bigr)=
\L_{|U^{\p\p\p}}\bigl(-(D_{|U^{\p\p}}+D^\p)\bigr).
$$
Now, let $\chi_*$ be the character of $H_{\ol{G}}(\la)$ as in Proposition~\ref{RamRamI}. By our strengthening of
the rationality properties, $\chi_*$ comes from a character of $H_{G}(\la)$ which we denote again by
$\chi_*$. We may view $\chi_*$ also as a character of $Q_{\GL(V)}(\la)$. 
The given $Q_{\GL(V)}(\la)$-bundle ${\cal Q}_{|U^{\p\p\p}}\subset {\cal P}_{|U^{\p\p\p}}$ and $\chi_*^{-1}$ define
a line bundle $\L_{\chi_*^{-1}}$, and ${\cal Q}_{|U^{\p\p\p}}$ and $W_{j_0}/W_{j_0-1}\otimes\C_{\chi_*^{-1}}$ define a vector bundle
$$
\widetilde{\cal F}^\vee\ \cong\ 
\bigl({\cal F}^\vee_{j_0|U^{\p\p\p}}/{\cal F}^\vee_{j_0-1|U^{\p\p\p}}\bigr)\otimes \L_{\chi_*^{-1}}
$$ over $U^{\p\p\p}$, so that 
$$
{\ol{\sigma}^{\p\p}}^*\bigl(\O_{\P(\widetilde{\cal F})}(1)\bigr)  \cong 
\L_{|U^{\p\p\p}}\bigl(-(D_{|U^{\p\p}}+D^\p)\bigr)\otimes \L^\vee_{\chi^{-1}_*}.
$$
Now, Proposition~\ref{Mummy} grants that the assumptions of Proposition~\ref{RamRamII} are satisfied,
so that we conclude that 
\begin{equation}
\label{BoundDegI}
 \deg(\L)\ge \deg\Bigl(\L_{|U^{\p\p\p}}\bigl(-(D_{|U^{\p\p}}+D^\p)\bigr)\Bigr) \ge 
\deg(\L_{\chi_*^{-1}}).
\end{equation}
By construction, there is a $q\in\Z_{>0}$ with $\chi_*^{-1}=q\cdot l_T(\la)$. 
Using Equation~(\ref{DegComp}), we find
$$
\deg(\L_{l_T(\la)}) = \sum_{i=1}^{s}\alpha_{i}\bigl(\deg(\E)\cdot \rk\E_i- \deg(\E_i)\cdot r\bigr).		
$$
The last expression is the coefficient of the monomial of degree ${\dim(X)-1}$ in $M(\E^\bullet,\ul{\alpha})$. 
Note that Remark \ref{FiniteInstabII} and Equation (\ref{BoundDegI}) imply that we
find a constant $C$, depending only on $a$, $b$, $c$, $d$, and $r$, such that
$$
C \ge \deg(\L_{l_T(\la)}).
$$
Let $\ol{\delta}$ be the coefficient of $x^{\dim X-1}$ in $\delta$ 
and observe $\mu(\E^\bullet,\ul{\alpha};\phi)
\le -1$. The semistability condition $M(\E^\bullet,\ul{\alpha})+\ol{\delta}\cdot\mu(\E^\bullet,\ul{\alpha};\phi)$ yields
the estimate
$$
C-\ol{\delta}\ge 0.
$$
Since this estimate was obtained assuming the negation of the assertion of the theorem,
our proof is finished here.
\end{proof}
\begin{Rem}
\label{MumfordSemStab0}
If $\deg(\E)=0=d$, then the above estimates show that we can choose $\delta_0=0$, i.e.,
whenever $\delta$ has degree exactly $\dim X-1$, we find $\widetilde{\sigma}_\eta\in
\P^{\rm ss}_{a,b,c}\times_\C \Spec(K)$.
\end{Rem}
\subsection{Step 2}
Now, let $\delta\in\Q[x]$ be a polynomial, such that the coefficient of $x^{\dim X-1}$ is larger
than $\delta_0$. We will call a sheaf $(\E,\phi)$ with a decoration of type $(a,b,c;\L)$ \it
asymptotically (semi)stable\rm, if a) 
$\widetilde{\sigma}_\eta\in \P^{\rm ss}_{a,b,c}\times_{\Spec(\C)} \Spec(K)$ 
(notation as in Theorem \ref{SemStabRed}) and b) for every weighted filtration
$(\E^\bullet, \ul{\alpha})$ with $\mu(\E^\bullet, \ul{\alpha};\phi)=0$, one has
$$
M(\E^\bullet, \ul{\alpha})(\succeq)0.
$$
The idea is, of course, that this notion of asymptotic (semi)stability equals the notion of
$\delta$-(semi)stability for ``large" $\delta$. The proof of this fact will be completed in Step 3.
First, we note a consequence of Theorem \ref{SemStabRed}.
\begin{Lem}
\label{EasyImp}
Under the above assumption on the parameter $\delta$, every $\delta$-(semi)\allowbreak stable sheaf $(\E,\phi)$
with a decoration of type $(a,b,c;\L)$ is asymptotically (semi)\allowbreak stable.
\end{Lem}
The next step is the following
\begin{Prop}
\label{StrongBoundedness}
For given non-negative integers $a$, $b$, $c$, a given line bundle $\L$ on $X$, and a given Hilbert polynomial 
$P$, the set of
isomorphy classes $[\E]$ of torsion free sheaves $\E$, such that there exists an asymptotically semistable
sheaf $(\E,\phi)$ with a decoration of type $(a,b,c;\L)$ is bounded.
\end{Prop}
\begin{proof}
Let ${\frak G}_i$ be the Gra\ss mannian of $i$-dimensional quotients of $\C^r$, $i=1,...,r-1$,
$r$ being the rank of an $\O_X$-module with Hilbert polynomial $P$.
Then, $\O_{{\frak G}_i}(1)$ stands for the polarization coming from the Pl\"ucker embedding
into $\P(W^i_1)$, $W^i_1:=\bigwedge^{r_i}\C^r$, $i=1,...,r-1$.
We fix integers $n_i\ge \eta_{i;\infty}$, $i=1,...,r-1$. Here, $\eta_{i;\infty}$ is the number
from Theorem \ref{InstFlagProduct} for the space $\P(W^i_1)\times\P(W_2)$, $W_2:=V_{a,b,c}$, $i=1,...,r-1$.
\par
Next, let $(\E,\phi)$ be an asymptotically semistable sheaf with a decoration of type
$(a,b,c;\L)$. Denote by $U$ the big open subset where $\E$ is locally free. Let $0\subsetneq\F\subsetneq \E$ be a saturated subsheaf. Then, there is a big open subset $U^\p\subseteq U$,
such that $\F_{|U^\p}$ is a subbundle of $\E_{|U^\p}$. Moreover, we may find a big open
subset $U^{\p\p}\subseteq U^\p$ and an effective divisor $D$ in $U^{\p\p}$, such that
$\phi_{|U^{\p\p}}$ induces a surjection 
$$
\phi^{\p\p}\colon \E_{a,b,c|U^{\p\p}}\lra \L_{|U^{\p\p}}(-D).
$$
Let ${\cal G}^{r-\rk\F}(\E_{|U^{\p\p}})\lra U^{\p\p}$ be the Gra\ss mann bundle of $(r-\rk\F)$-dimensional quotients of the fibres of $\E_{|U^{\p\p}}$. Note that the ample
line bundle $\O(1, n_{r-\rk\F})$ on ${\frak G}_{r-\rk\F}\times \P(V_{a,b,c})$ induces
a line bundle ${\cal H}_{r-\rk\F}$ on ${\cal G}^{r-\rk\F}(\E_{|U^{\p\p}})\times_{U^{\p\p}}
\P(\E_{a,b,c|U^{\p\p}})$. Furthermore, the subbundle $\F_{|U^{\p\p}}\subsetneq \E_{|U^{\p\p}}$
and $\phi^{\p\p}$ yield a section
$$
s^{\p\p}\colon U^{\p\p}\lra {\cal G}^{r-\rk\F}(\E_{|U^{\p\p}})\times_{U^{\p\p}}
\P(\E_{a,b,c|U^{\p\p}}),
$$ 
such that
$$
{s^{\p\p}}^*{\cal H}_{r-\rk\F}=\det(\E_{|U^{\p\p}}/\F_{|U^{\p\p}})\otimes \L_{|U^{\p\p}}(-D)^{\otimes n_{r-\rk\F}}.
$$
Therefore, if  $s^{\p\p}$ is generically semistable, we have
$$
\deg\bigl({s^{\p\p}}^*{\cal H}_{r-\rk\F}\bigr)\ge 0,
$$
by Proposition \ref{RamRamII}, i.e., 
\begin{equation}
\label{MaruyamaI}
\deg(\F)\le \deg(\E)+n_{r-\rk\F}\cdot\deg(\L)-n_{r-\rk\F}\cdot\deg(D)\le \deg(\E)+n_{r-\rk\F}\cdot\deg(\L).
\end{equation}
If the image of $s^{\p\p}$ is not generically semistable, then we set $W^\vee:= W^{r-\rk\F}_1\otimes S^{n_{r-\rk\F}}(W_2)$. As in the proof of Theorem \ref{SemStabRed}, we find
a one parameter subgroup $\widetilde{\la}\colon\C^*\lra \SL_r(\C)$, a subquotient $W^\p$ of
$W$ with the structure of an $L_{\GL_r(\C)}(\widetilde{\la})$-module structure, 
a big open subset $U^{\p\p\p}$, and a section
$$
s^{\p\p\p}\colon U^{\p\p\p}\lra \P\bigl(\E(W^\p)^\vee\bigr)
$$
with 
$$
{s^{\p\p\p}}^*({\cal H}_{W^\p})=\det(\E_{|U^{\p\p\p}}/\F_{|U^{\p\p\p}})\otimes \L_{|U^{\p\p\p}}^{\otimes n_{r-\rk\F}}
(-D^{\p\p}),
$$ 
for some
effective divisor $D^{\p\p}$ on $U^{\p\p\p}$ (composed of $n_{r-\rk\F}\cdot D_{|U^{\p\p\p}}$ and
some other effective divisor $D^\p$ on $U^{\p\p\p}$).
The vector bundle $\E(W^\p)$ corresponds to $\F^\vee_{j_0|U^{\p\p\p}}/\F_{j_0-1|U^{\p\p\p}}^\vee$ in the proof of Theorem \ref{SemStabRed},
and ${\cal H}_{W^\p}$ is the line bundle on $\P\bigl(\E(W^\p)^\vee\bigr)$ associated to
$\E(W^\p)$ and $\O_{\P(W^\p)}(1)$. There is again a character 
$\chi_*$ of $L_{\GL_r(\C)}(\widetilde{\la})$, such that $s^{\p\p\p}(\eta)$ is semistable
w.r.t.\ the linearization in ${\cal H}_{W^\p}$ twisted by the character $\chi_*$.
One concludes
\begin{equation}
\label{MaruyamaII}
\deg(\F)\le \deg(\E)+ n_{r-\rk\F}\cdot \deg(\L)-  \deg(D^{\p\p})+\deg(\L_{\chi_*}).
\end{equation}
Now, $\deg(\L_{\chi_*})$ is a \sl negative \rm multiple of the coefficient of the monomial of degree
$\dim X-1$ in $M(\E^\bullet,\ul{\alpha})$, the weighted filtration $(\E^\bullet,\ul{\alpha})$
being obtained as in the proof of Theorem \ref{SemStabRed} from an instability one parameter
subgroup. Theorem \ref{InstFlagProduct} implies
$$
\mu(\E^\bullet,\ul{\alpha};\phi)=0.
$$
The asymptotic semistability of $(\E,\phi)$, thus, grants $\deg(\L_{\chi_*})\le 0$, i.e.,
Inequality (\ref{MaruyamaII}) yields Inequality (\ref{MaruyamaI}).
We finally conclude by Maruyama's boundedness result, \cite{HL}, Theorem 3.3.7.
\end{proof}
\begin{Rem}
\label{MumfordSemStab}
If we are in the situation $\deg\E=0$ and $d=0$, the above arguments show that $\E$ must be a Mumford semistable
sheaf.
\end{Rem}
\subsection{Step 3}
\begin{Thm}
\label{Asymp}
There is a constant $\delta_1\ge \delta_0$, depending only on the input data $a$, $b$, $c$, $d$, and $P$,
such that for every stability parameter $\delta$ of degree exactly $\dim(X)-1$ the leading
coefficient of which is larger than $\delta_1$ and every 
torsion free sheaf $(\E,\phi)$ with a decoration of type $(a,b,c;\L)$ and with $P(\E)=P$, 
the following conditions are equivalent:
\begin{enumerate}
\item[\rm 1.] $(\E,\phi)$ is $\delta$-(semi)stable.
\item[\rm 2.] $(\E,\phi)$ is asymptotically (semi)stable.
\end{enumerate}
\end{Thm}
\begin{proof} The implication 1.\ $\Longrightarrow$ 2.\ has already been noted in Lemma \ref{EasyImp}.
For the converse, we use Proposition \ref{ConeDecompII}. This means that there is a finite set 
$$
{\cal T}=\bigl\{\,(r_1,...,r_{s_j};\alpha_1,...,\alpha_{s_j})\,|\, r_i\in \Z_{>0},\alpha_i\in\Q_{>0},
i=1,...,s_j,j=1,...,t\,\bigr\},
$$
such that the condition of $\delta$-(semi)stability has to be checked only for weighted filtrations
$(\E^\bullet,\ul{\alpha})$ with
\begin{equation}
\label{Stern}
(\rk \E_1,...,\rk\E_s;\alpha_1,...,\alpha_s)\in {\cal T}.
\end{equation}
Now, Equation (\ref{MaruyamaI}) and the finiteness of ${\cal T}$ imply that there is a constant $C$,
such that
$$
L(\E^\bullet,\ul{\alpha})\ge C
$$ 
for every asymptotically semistable decorated sheaf $(\E,\phi)$ (with the fixed input data)
and every weighted filtration $(\E^\bullet,
\ul{\alpha})$ of $\E$ which satisfies (\ref{Stern}).
We choose
$$
\delta_1\ge \max\{\, \delta_0, -C\,\}.
$$ 
\par
Suppose finally that $(\E,\phi)$ is an asymptotically (semi)stable torsion free sheaf with a decoration of type
$(a,b,c;\L)$ and with $P(\E)=P$ and that $\delta\in\Q[x]$ is a polynomial of degree exactly $\dim(X)-1$ the leading
coefficient $\ol{\delta}$ of which is larger than $\delta_1$. Recall that the definition of asymptotic semistability implies $\mu(\E^\bullet,\ul{\alpha};\phi)\ge 0$
for every asymptotically semistable decorated sheaf $(\E,\phi)$ and 
every weighted filtration $(\E^\bullet,\ul{\alpha})$ of $\E$. If $(\E^\bullet,\ul{\alpha})$ 
is a weighted filtration which satisfies 
(\ref{Stern}) and $\mu(\E^\bullet,\ul{\alpha};\phi)>0$,
then
$$
L(\E,\ul{\alpha})+\ol{\delta}\cdot\mu(\E^\bullet,\ul{\alpha};\phi)
> C+\delta_1\ge 0,
$$
so that $M(\E^\bullet,\ul{\alpha})+\delta\cdot \mu(\E^\bullet,\ul{\alpha};\phi)$ is indeed a positive polynomial.
If $(\E^\bullet,\ul{\alpha})$ is a weighted filtration which satisfies 
(\ref{Stern}) and $\mu(\E^\bullet,\ul{\alpha};\phi)=0$, then $M(\E^\bullet,\ul{\alpha})+\delta\cdot 
\mu(\E^\bullet,\ul{\alpha};\phi)$ 
is a positive (non-negative) polynomial, because $(\E,\phi)$ is asymptotically (semi)stable.
To conclude, we have verified the condition of $\delta$-(semi)stab\-il\-ity for all weighted filtrations $(\E^\bullet,
\ul{\alpha})$ of $\E$ which satisfy (\ref{Stern}), so that $(\E,\phi)$ is $\delta$-(semi)stable by Proposition
\ref{ConeDecompII}.
\end{proof}
\begin{Rem}
Suppose $\deg(\E)=0=d$.
By the Remarks \ref{MumfordSemStab0} and \ref{MumfordSemStab}, we may choose $\delta_0=0$ and $C=0$, whence
also $\delta_1=0$, so that, in this special case, Theorem \ref{Asymp} holds for every positive polynomial $\delta$
of degree exactly $\dim X-1$.
\end{Rem}
\begin{Cor}
Let $(\E,\phi)$ be a torsion free sheaf with a decoration of type $(a,b,c;\L)$, such that
$\widetilde{\sigma}_\eta$ lies in $\P^{\rm s}_{a,b,c}\times_{\Spec(\C)}\Spec(K)$
(notation as in {\rm Theorem \ref{SemStabRed}}), then $(\E,\phi)$ is asymptotically stable and,
in particular, $\delta$-stable for all $\delta\succ\hskip-3pt\succ 0$.
\end{Cor}
\begin{Rem}
This answers a question asked to me by Professor Okonek some time ago.
\end{Rem}
\begin{Ex}
Take $a=1$ and $c=0$, so that $(\E,\phi)$ can be viewed as a pair
consisting of a torsion free coherent sheaf $\E$ and a homomorphism $\widetilde{\phi}\colon \E\lra
\L^{\oplus b}$. Then, $\widetilde{\sigma}_\eta\in \P^{\rm s}_{a,b,c}\times_{\Spec(\C)}\Spec(K)$
is equivalent to the injectivity of $\widetilde{\phi}$. In that case the above theorem and corollary
were observed in \cite{OST}.
\end{Ex}
\subsection{Step 4}
Now, we carry out the final step in the proof of the main theorem. Let $\ol{\delta}_\infty>\delta_1$
with $\delta_1$ from Theorem \ref{Asymp}, and set $\delta_\infty :=\ol{\delta}_\infty\cdot x^{\dim X-1}$.
To begin with we note.
\begin{Prop}
\label{GlobalBound}
Fix the input data $a$, $b$, $c$, $\L$, and $P$. Then, the set ${\frak S}$ of isomorphism
classes $[\E]$ of torsion free sheaves $\E$ with Hilbert polynomial $P$
for which there exist a positive polynomial $\delta\in\Q[x]$
of degree at most $\dim X-1$ and a homomorphism $\phi\colon (\E^{\otimes a})^{\oplus b}\lra
\det(\E)^{\otimes c}\otimes \L$, such that the decorated sheaf $(\E,\phi)$ is $\delta$-semistable, is bounded.
\end{Prop}
\begin{proof}
If $(\E,\phi)$ is $\delta$-semistable for a polynomial $\delta\succeq\delta_\infty$, then
$(\E,\phi)$ is asymptotically semistable, and we have determined before a constant $C_1$
which depends only on the input data,
such that 
$$
\mu_{\rm max}(\E)\le \mu(\E)+C_1
$$
(see Proposition \ref{StrongBoundedness}).
If, on the other hand, $\delta\preceq\delta_\infty$, then, letting $\ol{\delta}$ be the coefficient of
$x^{\dim X-1}$ in $\delta$, we find
$$
\mu_{\rm max}(\E)\le \mu(\E)+\frac{\ol{\delta}\cdot a\cdot (r-1)}{r}\le 
\mu(\E)+\frac{\ol{\delta}_\infty\cdot a\cdot (r-1)}{r}.
$$
This is proved
in \cite{Schmitt0}, Theorem 2.6.
Since the constant $C_2:=\max\{\, 0,C_1, (\ol{\delta}\cdot a\cdot (r-1))/r\,\}$ depends only on the input data, 
the assertion is again a consequence of Maruyama's boundedness theorem (\cite{HL}, Theorem 3.3.7).
\end{proof}
Note that we may rewrite the condition on $\mu_{\max}(\E)$ we have just found as
\begin{equation}
\label{LowerBound}
\deg(\E)\cdot \rk\F-\deg(\F)\cdot\rk\F\ge -C_2\cdot\rk\F\cdot\rk\E\ge -C_2^\p:=-C_2\cdot (r-1)\cdot r.
\end{equation}
Let ${\cal T}$ be the set defined in the discussion before Proposition \ref{ConeDecompII}. Next, we note
that for a weighted filtration $(\E^\bullet,\ul{\alpha})$ of $\E$, one
has
\begin{equation}
\label{LowerBoundII}
\mu(\E^\bullet, \ul{\alpha};\phi)\ge -\sum_{i=1}^s\alpha_i(r-1),
\end{equation}
by \cite{Schmitt0}, Lemma 1.8 (ii).
Using (\ref{LowerBound}) and (\ref{LowerBoundII}), it is easy to determine a constant $C_3>0$,
such that for every torsion free sheaf $\E$ and every weighted filtration $(\E^\bullet,\ul{\alpha})$,
such that
\begin{itemize}
\item $\E\in {\frak S}$,
\item $(\rk\E_1,...,\rk\E_s,\alpha_1,...,\alpha_s)\in {\cal T}$, and
\item there exists an index $i_0\in \{\, 1,...,s\,\}$ with $\mu(\E_{i_0})<\mu(\E)-C_3$,
\end{itemize}
one has
\begin{equation}
\label{BoundII}
L(\E^\bullet,\ul{\alpha})>\ol{\delta}_\infty\cdot \sum_{i=1}^s\alpha_i(r-1).
\end{equation}
On the other hand, the set ${\frak S}(C_3)$ of isomorphism classes $[\F]$ of subsheaves $0\subsetneq\F\subsetneq\E$
with $[\E]\in {\frak S}$ and $\mu(\F)\ge \mu(\E)-C_3$ is bounded, too. Let
$$
{\cal U}:=\bigl\{\,P_1,...,P_u\,\bigr\}
$$ 
be the set of polynomials of the form $M(\E^\bullet,\ul{\alpha})$ where
$(\E^\bullet, \ul{\alpha})$ is a weighted filtration of a torsion free sheaf $\E$, such that
\begin{itemize}
\item $\E\in {\frak S}$,
\item $(\rk\E_1,...,\rk\E_s,\alpha_1,...,\alpha_s)\in {\cal T}$, and
\item $\mu(\E_{i_0})\ge\mu(\E)-C_3$, $i=1,...,s$.
\end{itemize}
Now, (\ref{BoundII}) yields the following
\begin{Lem}
\label{BoundIII}
Suppose we are given a torsion free sheaf $(\E,\phi)$, $[\E]\in {\frak S}$, with a decoration of type $(a,b,c;\L)$,
a polynomial $\delta\in \Q[x]$ with $0\prec\delta\preceq\ol{\delta}_\infty$, and a weighted
filtration $(\E^\bullet,\ul{\alpha})$ with $(\rk\E_1,...,\rk\E_s,\alpha_1,...,\alpha_s)\in {\cal T}$, such that
$$
M(\E^\bullet,\ul{\alpha})+\delta\cdot \mu(\E^\bullet, \ul{\alpha};\phi)\prec 0.
$$
Then,
$$
M(\E^\bullet,\ul{\alpha})\in {\cal U}.
$$
\end{Lem}
Next, let
$$
{\cal V}:=\bigl\{\, z_1,...,z_v\,\bigr\}\subset\Z
$$
be the set of finite values for $\mu(\la, w)$ where
\begin{itemize}
\item $w\in \P(V_{a,b,c})$,
\item $(\dim V_1,...,\dim V_s, \alpha_1,...,\alpha_s)\in {\cal T}$, $(V^\bullet,\ul{\alpha})$ being the weighted
flag of $\la$.
\end{itemize}
Finally, suppose that $\delta_1\preceq\delta_\infty$ and $\delta_2\preceq\delta_\infty$ 
are two positive polynomials and that $(\E,\phi)$ is a torsion free sheaf with a decoration of type
$(a,b,c;\L)$, such that $P(\E)=P$, and suppose that $(\E,\phi)$ is $\delta_1$-(semi)stable
(whence $[\E]\in {\frak S}$), but 
not $\delta_2$-(semi)stable. By Proposition \ref{ConeDecompII}, there must be a weighted filtration
$(\E^\bullet, \ul{\alpha})$ with $(\rk\E_1,...,\rk\E_s,\allowbreak \alpha_1,...,\alpha_s)\in {\cal T}$, such that
\begin{eqnarray*}
M(\E^\bullet,\ul{\alpha})+\delta_1\cdot \mu(\E^\bullet,\ul{\alpha};\phi) &\succeq (\succ)& 0
\\
M(\E^\bullet,\ul{\alpha})+\delta_2\cdot \mu(\E^\bullet,\ul{\alpha};\phi) &\prec (\preceq)& 0.
\end{eqnarray*}
Suppose first that $\delta_1\prec\delta_2$. Then, necessarily $\mu(\E^\bullet,\ul{\alpha};\phi)<0$
and $M(\E^\bullet,\ul{\alpha})\succ 0$. We see
\begin{equation}
\label{BoundIV}
\delta_1\preceq-\frac{M(\E^\bullet,\ul{\alpha})}{\mu(\E^\bullet,\ul{\alpha};\phi)}\preceq\delta_2.
\end{equation}
If, on the other hand, $\delta_2\prec\delta_1$, we must have $\mu(\E^\bullet,\ul{\alpha};\phi)>0$
and $M(\E^\bullet,\ul{\alpha})\prec 0$, and then
\begin{equation}
\label{BoundV}
\delta_1\succeq-\frac{M(\E^\bullet,\ul{\alpha})}{\mu(\E^\bullet,\ul{\alpha};\phi)}\succeq\delta_2.
\end{equation}
Let
$$
\left\{\, -\frac{Q}{z}\succ 0\,\Big|\, Q\in {\cal U}, z\in {\cal V}\,\right\}
=\bigl\{\,\widetilde{\delta}_1,...,\widetilde{\delta}_w\,\bigr\}
\q\hbox{with}\ \widetilde{\delta}_1\prec\cdots\prec\widetilde{\delta}_w.
$$
To conclude, we may summarize (\ref{BoundIV}), (\ref{BoundV}), Lemma \ref{BoundIII}, and some
easy thoughts in
\begin{Lem} 
Let
$(\E,\phi)$ be a torsion free sheaf 
with a decoration of type $(a,b,c;\allowbreak \L)$, such that $P(\E)=P$.
\par
{\rm i)} Suppose there is an index $i_0\in \{\, 0,...,w-1\,\}$, such that $\widetilde{\delta}_{i_0}
\prec\delta_1\prec\delta_2\prec\widetilde{\delta}_{i_0+1}$, then $(\E,\phi)$
is $\delta_1$-(semi)stable if and only if
it is $\delta_2$-(semi)stable. Here, we have set $\widetilde{\delta}_0:=0$.
\par
{\rm ii)} If $\delta\succ\widetilde{\delta}_w$, then $(\E,\phi)$ is $\delta$-(semi)stable, if and only if
it is asymptotically (semi)stable.
\par
{\rm iii)} Suppose $\delta\prec\delta_1$. If $(\E,\phi)$ is $\delta$-semistable, then $\E$ is a semistable sheaf. 
If $\E$ is stable sheaf, then $(\E,\phi)$ is $\delta$-stable.
\par
{\rm iv)} Suppose $\widetilde{\delta}_i\prec\delta\prec\widetilde{\delta}_{i+1}$, for some $i\in \{\,0,...,w\,\}$.
If $(\E,\phi)$ is $\delta$-semistable, then $(\E,\phi)$ is also $\widetilde{\delta}_i$- and 
$\widetilde{\delta}_{i+1}$-semistable.
If $(\E,\phi)$ is $\widetilde{\delta}_i$-stable or $\widetilde{\delta}_{i+1}$-stable, 
then it is likewise $\delta$-stable.
\end{Lem}
By this lemma, the main theorem holds for a finite subset $\{\,\widehat{\delta}_1,...,\widehat{\delta}_m\,\}$
of $\{\,\widetilde{\delta}_1,...,\widetilde{\delta}_w\,\}$.\qed

\end{document}